\documentclass[a4paper]{amsart}

\usepackage[mathscr]{eucal}
\usepackage{amsmath,amsfonts,amssymb}
\usepackage{amsthm}
\usepackage[dvips]{graphicx,color}
\graphicspath{{./figure/}}
\usepackage{bm}

\newcommand{\mathsym}[1]{{}}

\theoremstyle{plain}
\newtheorem{theorem0}{Theorem}

\newtheorem{proposition0}[theorem0]{Proposition}
\newtheorem{lemma0}[theorem0]{Lemma}
\newtheorem{theorem}{Theorem}[section]
\newtheorem{proposition}[theorem]{Proposition}

\theoremstyle{definition}
\newtheorem{definition0}[theorem0]{Definition}

\begin{document}

\title[Certain lattice neighboring from Barnes-Wall lattice and some laminated lattices]{Certain lattice neighboring\\ from Barnes-Wall lattice and some laminated lattices}
\author{Junichi Shigezumi}
\date{}

\maketitle \vspace{-0.1in}
\begin{center}
Graduate School of Mathematics Kyushu University\\
Hakozaki 6-10-1 Higashi-ku, Fukuoka, 812-8581 Japan\\
{\it E-mail address} : j.shigezumi@math.kyushu-u.ac.jp \vspace{-0.05in}
\end{center} \quad

\begin{quote}
{\small\bfseries Abstract.}
We give classifications of integral lattices which include the Barnes-Wall lattice $BW_{16}$ or laminated lattices of dimension $1 \leqslant d \leqslant 8$ and of minimum $4$. Also, we give certain lattice neighboring from each lattice. Furthermore, we study spherical designs and the other properties.\\  \vspace{-0.15in}

\noindent
{\small\bfseries Key Words and Phrases.}
lattice neighboring, Euclidean lattice, spherical design.\\ \vspace{-0.15in}

\noindent
2000 {\it Mathematics Subject Classification}. Primary 05B30; Secondary 03G10. \vspace{0.15in}
\end{quote}

\section*{Introduction}

We have the following {\it lattice neighboring} from the Barnes-Wall lattice $BW_{16}$ to the three unimodular lattices $\Lambda_{16, 2, 3}$, $\bm{D}_{16}^{+}$, and $\bm{E}_{8} \perp \bm{E}_{8}$:\vspace{-0.25in}

\begin{figure}[h]
\includegraphics[width=4.5in]{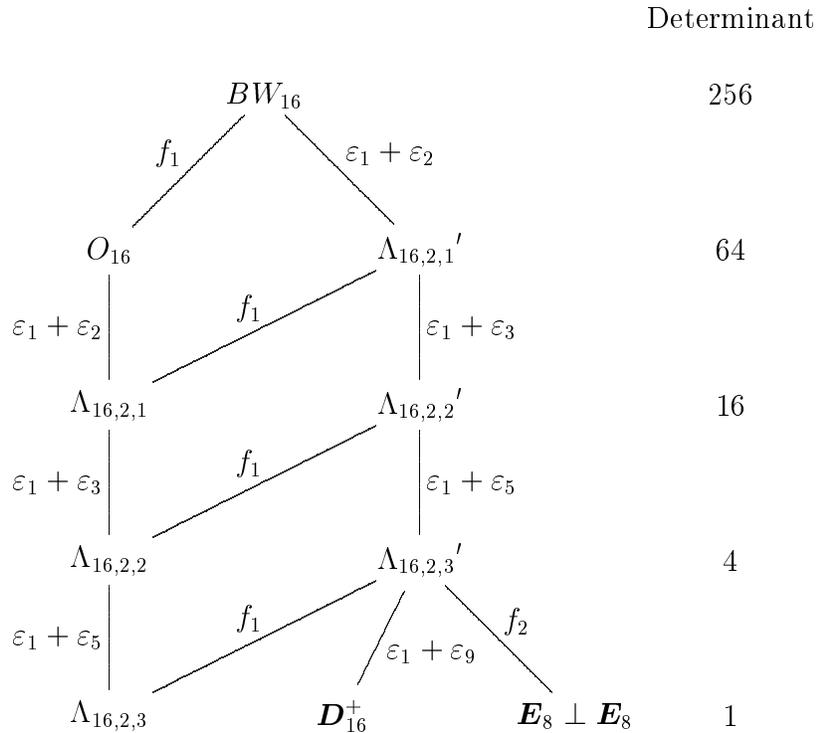}\vspace{-0.25in}
\caption{Lattice neighboring from $BW_{16}$}\label{fig16}
\end{figure}

Here, the above figure means that $O_{16} =  BW_{16} \cup (f_1 + BW_{16})$, for example. We will introduce the definitions of the lattices and vectors for the lattice neighboring in Figure $\ref{fig16}$ are in Section \ref{sec-def16}.

{\it Lattice neighboring} is to consider {\it neighbors} of the given lattice, whose intersection has index two in each of them. In this note, we consider lattice neighborings with integral lattices, that is we consider an integral lattice $L' = \langle L, x\rangle$ as a neighbor of the given lattice $L$ and some vector $x \in L^{\sharp} \setminus L$ of integral norm if $\langle L, x\rangle = L \cup (x + L)$.

\makeatletter
\renewcommand{\thetheorem}{1.1}
\makeatletter

Furthermore, we have the following fact:
\begin{theorem}[see Section 1]
The integral lattices which include $BW_{16}$ are isometric to one of the nine lattices in Figure $\ref{fig16}$, namely $BW_{16}$, $O_{16}$, $\Lambda_{16, 2, 1}$, $\Lambda_{16, 2, 2}$, $\Lambda_{16, 2, 3}$, ${\Lambda_{16, 2, 1}}'$, ${\Lambda_{16, 2, 2}}'$, ${\Lambda_{16, 2, 3}}'$, $\bm{D}_{16}^{+}$, and $\bm{E}_{8} \perp \bm{E}_{8}$.
\end{theorem}\quad

\makeatletter
\renewcommand{\thetheorem}{\arabic{section}.\arabic{theorem}}
\@addtoreset{theorem}{section}
\makeatletter

Similarly, we consider lattice neighborings for the following {\it laminated lattices}.

We denote by $\Lambda_n$ the laminated lattices of minimum $4$ for $0 \leq n \leq 24$. We can refer to \cite{CS} for the definition of laminated lattices. We have the following theorem:
\begin{proposition0}[cf. \cite{CS} Ch.5, Theorem 2]
The laminated lattices $\Lambda_0, \Lambda_1, \ldots, \Lambda_8$ are unique and are isomorphic to $\Lambda_0, \bm{A}_1 \cong \mathbb{Z}, \bm{A}_2, \bm{A}_3 \cong \bm{D}_3, \bm{D}_4, \bm{D}_5, \bm{E}_6, \bm{E}_7, \bm{E}_8$ respectively. Their determinants are shown in the following table:\\

\begin{center}
\begin{tabular}{r|ccccccccc}
\hline
$n$ & $0$ & $1$ & $2$ & $3$ & $4$ & $5$ & $6$ & $7$ & $8$\\
\hline \hline
Lattice & $\Lambda_0$ & $2 \, \mathbb{Z}$ & $\sqrt{2} \, \bm{A}_2$ & $\sqrt{2} \, \bm{A}_3$ & $\sqrt{2} \, \bm{D}_4$ & $\sqrt{2} \, \bm{D}_5$ & $\sqrt{2} \, \bm{E}_6$ & $\sqrt{2} \, \bm{E}_7$ & $\sqrt{2} \, \bm{E}_8$\\
\hline
Det. & $1$ & $4$ & $12$ & $32$ & $64$ & $128$ & $192$ & $256$ & $256$\\
\hline
\end{tabular}
\end{center}\quad
\end{proposition0}
\noindent
Furthermore, $\Lambda_{16}$ is isometric to $BW_{16}$, and that $\Lambda_{24}$ is the {\it Leech lattice}.

If the determinant of $\Lambda_d$ is a power of four, then, similarly to $BW_{16}$, we can consider a lattice neighboring from $\Lambda_d$ to some unimodular lattices.\\

Let $S^{d-1} := \{ (x_1, \ldots, x_d) \in \mathbb{R}^d \, ; \, x_1^2 + \cdots + x_d^2 = 1 \}$ be the Euclidean sphere for $d \geqslant 1$.

\begin{definition0}[Spherical design \cite{DGS}]\label{def-design}
Let $X$ be a non-empty finite set on the Euclidean sphere $S^{d-1}$, and let $t$ be a positive integer. $X$ is called a spherical $t$-design if
\begin{equation}
\frac{1}{| S^{d-1} |} \int_{S^{d-1}} f(\xi) \, d \xi \ = \ \frac{1}{| X |} \hspace{0.05in} \sum_{\xi \in X} \hspace{0.05in} f(\xi) \label{eq-design}
\end{equation}
for every polynomial $f(x) = f(x_1, \ldots, x_d)$ of degree at most $t$.
\end{definition0}

Here, the left-hand side of the above equation (\ref{eq-design}) means the average on the sphere $S^{d-1}$, and the right-hand side means the average on the finite subset $X$. Thus, if $X$ is a spherical design, then $X$ gives a certain approximation of the sphere $S^{d-1}$.

Furthermore, we denote the distance set of $X$ by $A(X) := \{ (x, y) \: ; \: x, y \in X, x \ne y \}$; then we call $X$ an $s$-distance set if $| A(X) | = s$. Now, $X$ is said to be a $(d, n, s, t)$-configuration if $X \subset S^{d-1}$ is of order $n (:= | X |)$, a $s$-distance set, and a spherical $t$-design.\\

The set $s_m (L)$ of vectors of the lattice $L$ which take the same value $m$ for their norm is called the {\it shell of the lattice}, i.e. $s_m (L) := \{ x \in L \: ; \: (x, x) = m \}$. For every nonempty shell $s_m (L)$ of a lattice $L$, a normalization $X = \frac{1}{\sqrt{m}} s_m (L)$ is considered, where $X$ is a finite set on an Euclidean sphere. A lattice, whose minimal shell is a spherical $4$-design (i.e. a $5$-design), is said to be {\it strongly perfect}.

\setcounter{section}{-1}
\section{Preliminaries}\label{sec-def0}
In this note, we denote by $\varepsilon_1, \ldots, \varepsilon_d$ an orthonormal basis of $\mathbb{R}^d$, and we denote by $e_1, \ldots, e_d$ a basis of the lattice $L = \mathbb{Z} e_1 \oplus \cdots \oplus \mathbb{Z} e_d \subset \mathbb{R}^d$. Furthermore, we also denote $L^m := L \perp L \perp \cdots \perp L$ and $L(x_1, \ldots, x_k) := \langle L, x_1, \ldots, x_k \rangle$.

B. B. Venkov introduces the following construction of the lattice:

\begin{lemma0}[Venkov \cite{V}, Lemma 7.1]
Let $L$ be an even integral lattice of dimension $n \geq 2$ and of minimum $4$, and let $e$ be a minimal vector of $L$. Denote by $p$ the orthogonal projection on the hyperplane $H = e^{\perp}$, put ${L_e}' = \{ x \in L \ | \ (e, x) \equiv 0 \pmod{2} \}$, and let $L_e = p({L_e}')$. Suppose that one of the following two assumptions holds:
\def\labelenumi{(\arabic{enumi})}
\begin{enumerate}
\item There is $x \in L$ such that $(e, x) \equiv 1 \pmod{2}$;

\item We have $(y, e) \equiv 0 \pmod{2}$ for all $y \in L$, and $L$ contains a vector $x$ such that $(e, x) \equiv 2 \pmod{4}$.
\end{enumerate}
Then, $L_e$ is a odd integral lattice of minimum at least $3$, and we have $\det(L_e) = \det(L)$ under assumption $(1)$ and $\det(L_e) = \frac{1}{4} \det(L)$ under assumption $(2)$.
\end{lemma0}

Now, we set $O_1 = \sqrt{3} \mathbb{Z}$. We denote by $O_7$ (resp. $O_{23}$) the projected $L_e$ associated with the laminated lattice $\Lambda_8$ (resp. $\Lambda_{24}$). Finally, we denote by $O_{22}$ (resp. $O_{16}$) the orthogonal of $O_1$ (resp. $O_7$) in $O_{23}$. Only these five lattice are the strongly perfect lattices of minimum $3$, classified by B. B. Venkov \cite{V}.\\

For an orthonormal basis $\varepsilon_1, \ldots, \varepsilon_d$ of $\mathbb{R}^d$, we often denote
\begin{align*}
(A_1)^{2 k} &:= \{ \pm (\varepsilon_{2 i - 1} \pm \varepsilon_{2 i}) \: ; \: 1 \leqslant i \leqslant k \}
 & (\bm{A}_1)^{2 k} & := \langle (A_1)^{2 k} \rangle\\
D_d &:= \{ \pm (\varepsilon_i \pm \varepsilon_j) \: ; \: 1 \leqslant i < j \leqslant d \}
 & \bm{D}_d &:= \langle D_d \rangle\\
(D_4)^{k} &:= \{ \pm (\varepsilon_i \pm \varepsilon_j) \: ; \: 4 k' - 3 \leqslant i < j \leqslant 4 k' \; \text{for} \; 1 \leqslant k' \leqslant k \}
 & (\bm{D}_4)^{k} &:= \langle (D_4)^{k} \rangle\\
(D_8)^{k} &:= \{ \pm (\varepsilon_i \pm \varepsilon_j) \: ; \: 8 k' - 7 \leqslant i < j \leqslant 8 k' \; \text{for} \; 1 \leqslant k' \leqslant k \}
 & (\bm{D}_8)^{k} &:= \langle (D_8)^{k} \rangle\\
f & := \frac{\varepsilon_1 + \cdots + \varepsilon_8}{2}
 & \bm{E}_8 &:= \langle D_8, f \rangle\\
\end{align*}\quad

\section{the Barnes-Wall lattice $BW_{16}$}

We will introduce the definitions of the lattices and vectors for the lattice neighboring in Figure $\ref{fig16}$ are in the next section.

It is known that, in $16$-dimensional space $\mathbb{R}^{16}$, there are just three unimodular lattices of minimum at least $2$ (cf. \cite{CS} Ch. 16).  Furthermore, we can show the following fact:
\begin{theorem}\label{prop-intlat16}
The integral lattices which include $BW_{16}$ are isometric to one of the nine lattices in Figure $\ref{fig16}$, namely $BW_{16}$, $O_{16}$, $\Lambda_{16, 2, 1}$, $\Lambda_{16, 2, 2}$, $\Lambda_{16, 2, 3}$, ${\Lambda_{16, 2, 1}}'$, ${\Lambda_{16, 2, 2}}'$, ${\Lambda_{16, 2, 3}}'$, $\bm{D}_{16}^{+}$, and $\bm{E}_{8} \perp \bm{E}_{8}$.
\end{theorem}

It is easy to classify the above lattices. It is well known that the Barnes-Wall lattice $BW_{16}$ is {\it $2$-elementary}, where we have $2 BW_{16}^{\sharp} \subset BW_{16}$. Furthermore, we have the following:

\begin{proposition}
The $2^8$ classes of $BW_{16}^{\sharp} / BW_{16}$ consisit of:
\begin{trivlist}
\item[$(i)$] the zero class,

\item[$(ii)$] $135$ classes each represented by a vector of $BW_{16}^{\sharp}$ of norm $2$,

\item[$(ii)$] $120$ classes each represented by a vector of $BW_{16}^{\sharp}$ of norm $3$.
\end{trivlist}
\end{proposition}

Thus, to classify such lattices, we have only to consider the neighboring with the each vector of norm $2$ and $3$. We can easily check by computer search. In fact, we have $19381$ integral lattices from the combinations of the above classes. The following table contains the number of integral lattices which are isometric to each lattice in Proposition \ref{prop-intlat16}:\\

\begin{center}
\begin{tabular}{c|ccc|ccc|ccc|ccc|c}
\cline{1-2} \cline{4-5} \cline{7-8} \cline{10-11} \cline{13-14}
$BW_{16}$ & $1$ & \qquad & $O_{16}$ & $120$ & \qquad & $\Lambda_{16, 2, 1}$ & $3780$ & \qquad & $\Lambda_{16, 2, 2}$ & $9450$ & \qquad & $\Lambda_{16, 2, 3}$ & $2025$\\
  &&& ${\Lambda_{16, 2, 1}}'$ & $135$ && ${\Lambda_{16, 2, 2}}'$ & $1575$ && ${\Lambda_{16, 2, 3}}'$ & $2025$ && $\bm{D}_{16}^{+}$ & $135$\\
  &&&&&&&&&&&& $\bm{E}_{8} \perp \bm{E}_{8}$ & $135$\\
\cline{1-2} \cline{4-5} \cline{7-8} \cline{10-11} \cline{13-14}
\end{tabular}
\end{center}\quad

In the following sections, we will introduce the definitions and some properties of lattices.\\

\subsection{Definitions}\label{sec-def16}

Here, we denote some vectors
\begin{gather*}
f_0 := \frac{- \varepsilon_1 + \varepsilon_2 + \cdots + \varepsilon_8 - \varepsilon_9 + \varepsilon_{10} + \cdots + \varepsilon_{16}}{2},\\
f_1 := \frac{\varepsilon_1 + \cdots + \varepsilon_8}{2} + \varepsilon_9, \quad
f_2 := \frac{- \varepsilon_1 + \cdots + \varepsilon_8}{2},
\end{gather*}
\begin{align*}
g_1 &:= \varepsilon_1 + \varepsilon_5 + \varepsilon_9 + \varepsilon_{13}, &
g_2 &:= \varepsilon_1 + \varepsilon_3 + \varepsilon_5 + \varepsilon_7, &
g_3 &:= \varepsilon_1 + \varepsilon_3 + \varepsilon_9 + \varepsilon_{11}, \\
g_4 &:= \varepsilon_1 + \varepsilon_3 + \varepsilon_{13} + \varepsilon_{15}, &
g_5 &:= \varepsilon_1 + \varepsilon_2 + \varepsilon_3 + \varepsilon_4, &
g_6 &:= \varepsilon_1 + \varepsilon_2 + \varepsilon_5 + \varepsilon_6, \\
g_7 &:= \varepsilon_1 + \varepsilon_2 + \varepsilon_7 + \varepsilon_8,  &
g_8 &:= \varepsilon_1 + \varepsilon_2 + \varepsilon_9 + \varepsilon_{10}, &
g_9 &:= \varepsilon_1 + \varepsilon_2 + \varepsilon_{11} + \varepsilon_{12}, \\
g_{10} &:= \varepsilon_1 + \varepsilon_2 + \varepsilon_{13} + \varepsilon_{14},  &
g_{11} &:= \varepsilon_1 + \varepsilon_2 + \varepsilon_{15} + \varepsilon_{16}.
\end{align*}\quad

First, we have the following definition for the Barnes-Wall lattice $BW_{16}$:
\begin{equation}
BW_{16} := \langle 2 \varepsilon_1, 2 \varepsilon_2, \ldots, 2 \varepsilon_{16}, f_0, g_1, \ldots, g_{11} \rangle
\end{equation}
We also have the following {\it generator matrix}:
{\footnotesize \begin{equation}
\frac{1}{2}
\left(\begin{array}{cccccccccccccccc}
 4 & 0 & 0 & 0 & 0 & 0 & 0 & 0 & 0 & 0 & 0 & 0 & 0 & 0 & 0 & 0 \\
 0 & 4 & 0 & 0 & 0 & 0 & 0 & 0 & 0 & 0 & 0 & 0 & 0 & 0 & 0 & 0 \\
 0 & 0 & 4 & 0 & 0 & 0 & 0 & 0 & 0 & 0 & 0 & 0 & 0 & 0 & 0 & 0 \\
 2 & 2 & 2 & 2 & 0 & 0 & 0 & 0 & 0 & 0 & 0 & 0 & 0 & 0 & 0 & 0 \\
 2 & 2 & 0 & 0 & 2 & 2 & 0 & 0 & 0 & 0 & 0 & 0 & 0 & 0 & 0 & 0 \\
 0 & 0 & 0 & 0 & 0 & 4 & 0 & 0 & 0 & 0 & 0 & 0 & 0 & 0 & 0 & 0 \\
 2 & 0 & 2 & 0 & 2 & 0 & 2 & 0 & 0 & 0 & 0 & 0 & 0 & 0 & 0 & 0 \\
 2 & 2 & 0 & 0 & 0 & 0 & 2 & 2 & 0 & 0 & 0 & 0 & 0 & 0 & 0 & 0 \\
 2 & 2 & 0 & 0 & 0 & 0 & 0 & 0 & 2 & 2 & 0 & 0 & 0 & 0 & 0 & 0 \\
 2 & 0 & 2 & 0 & 0 & 0 & 0 & 0 & 2 & 0 & 2 & 0 & 0 & 0 & 0 & 0 \\
 2 & 2 & 0 & 0 & 0 & 0 & 0 & 0 & 0 & 0 & 2 & 2 & 0 & 0 & 0 & 0 \\
 2 & 0 & 0 & 0 & 2 & 0 & 0 & 0 & 2 & 0 & 0 & 0 & 2 & 0 & 0 & 0 \\
 2 & 2 & 0 & 0 & 0 & 0 & 0 & 0 & 0 & 0 & 0 & 0 & 2 & 2 & 0 & 0 \\
 2 & 0 & 2 & 0 & 0 & 0 & 0 & 0 & 0 & 0 & 0 & 0 & 2 & 0 & 2 & 0 \\
 2 & 2 & 0 & 0 & 0 & 0 & 0 & 0 & 0 & 0 & 0 & 0 & 0 & 0 & 2 & 2 \\
 -1 & 1 & 1 & 1 & 1 & 1 & 1 & 1 & -1 & 1 & 1 & 1 & 1 & 1 & 1 & 1
\end{array}\right)
\end{equation}
}Incidentally, we have the following more beautiful generator matrix of $BW_{16}$ with another construction (cf. \cite{CS} Ch. 4):
{\scriptsize \begin{equation}
\frac{1}{\sqrt{2}}
\left(\begin{array}{cccccccccccccccc}
 4 & 0 & 0 & 0 & 0 & 0 & 0 & 0 & 0 & 0 & 0 & 0 & 0 & 0 & 0 & 0 \\
 2 & 2 & 0 & 0 & 0 & 0 & 0 & 0 & 0 & 0 & 0 & 0 & 0 & 0 & 0 & 0 \\
 2 & 0 & 2 & 0 & 0 & 0 & 0 & 0 & 0 & 0 & 0 & 0 & 0 & 0 & 0 & 0 \\
 2 & 0 & 0 & 2 & 0 & 0 & 0 & 0 & 0 & 0 & 0 & 0 & 0 & 0 & 0 & 0 \\
 2 & 0 & 0 & 0 & 2 & 0 & 0 & 0 & 0 & 0 & 0 & 0 & 0 & 0 & 0 & 0 \\
 2 & 0 & 0 & 0 & 0 & 2 & 0 & 0 & 0 & 0 & 0 & 0 & 0 & 0 & 0 & 0 \\
 2 & 0 & 0 & 0 & 0 & 0 & 2 & 0 & 0 & 0 & 0 & 0 & 0 & 0 & 0 & 0 \\
 2 & 0 & 0 & 0 & 0 & 0 & 0 & 2 & 0 & 0 & 0 & 0 & 0 & 0 & 0 & 0 \\
 2 & 0 & 0 & 0 & 0 & 0 & 0 & 0 & 2 & 0 & 0 & 0 & 0 & 0 & 0 & 0 \\
 2 & 0 & 0 & 0 & 0 & 0 & 0 & 0 & 0 & 2 & 0 & 0 & 0 & 0 & 0 & 0 \\
 2 & 0 & 0 & 0 & 0 & 0 & 0 & 0 & 0 & 0 & 2 & 0 & 0 & 0 & 0 & 0 \\
 1 & 1 & 1 & 1 & 0 & 1 & 0 & 1 & 1 & 0 & 0 & 1 & 0 & 0 & 0 & 0 \\
 0 & 1 & 1 & 1 & 1 & 0 & 1 & 0 & 1 & 1 & 0 & 0 & 1 & 0 & 0 & 0 \\
 0 & 0 & 1 & 1 & 1 & 1 & 0 & 1 & 0 & 1 & 1 & 0 & 0 & 1 & 0 & 0 \\
 0 & 0 & 0 & 1 & 1 & 1 & 1 & 0 & 1 & 0 & 1 & 1 & 0 & 0 & 1 & 0 \\
 1 & 1 & 1 & 1 & 1 & 1 & 1 & 1 & 1 & 1 & 1 & 1 & 1 & 1 & 1 & 1
\end{array}\right)
\end{equation}}

Secondly, we have another definition $O_{16} := \langle 2 \varepsilon_1, 2 \varepsilon_2, \ldots, 2 \varepsilon_{16}, f_0, f_1, g_1, \ldots, g_{11} \rangle$.

Thirdly, we define the three lattices, namely $\Lambda_{16,2,1} := \langle (A_1)^{16}, f_0, f_1, g_1, g_2, g_3, g_4 \rangle$, $\Lambda_{16,2,2} := \langle (D_4)^4, f_0, f_1, g_1 \rangle$, and $\Lambda_{16,2,3} := \langle D_8 \perp D_8, f_0, f_1 \rangle$, where root systems in the above definitions are defined as Section \ref{sec-def0}. Then, we have $\langle 2 \varepsilon_1, \ldots, 2 \varepsilon_{16} \rangle \subset (A_1)^{16} \subset (D_4)^4 \subset (D_8)^2$. These lattices are of dimension $16$ and of minimum $2$, and its shell of norm $3$ is spherical $5$ design. Such lattices are classified in \cite{S}, there exist just three.

Similarly, we define the three lattices, namely ${\Lambda_{16,2,1}}' := \langle (A_1)^{16}, f_0, g_1, g_2, g_3, g_4 \rangle$, ${\Lambda_{16,2,2}}' := \langle (D_4)^4, f_0, g_1 \rangle$, and ${\Lambda_{16,2,3}}' := \langle D_8 \perp D_8, f_0 \rangle$.

Finally, we define the two even unimodular lattice $\bm{D}_{16}^{+} = \ \langle D_{16}, f_0 \rangle \ := \langle D_8 \perp D_8, \varepsilon_1 + \varepsilon_9, f_0 \rangle$ and $\bm{E}_{8} \perp \bm{E}_{8} = \langle E_8 \perp E_8 \rangle := \langle D_8 \perp D_8, f_0, f_1 \rangle$.\\

\subsection{Properties of lattices}

The theta series of each lattice have the following form:
\begin{align*}
\Theta_{BW_{16}} &= 1 + 4320 \, q^4 + 61440 \, q^6 + 522720 \, q^8 + 2211840 \, q^{10} + 8960640 \, q^{12} + \cdots\\
\quad\\
\Theta_{O_{16}} &= 1 \hspace{0.45in} + 512 \, q^3 + 4320 \, q^4 + 18432 \, q^5 + 61440 \, q^6 + 193536 \, q^7\\
 &+ 522720 \, q^8 + 1126400 \, q^9 + 2211840 \, q^{10} + 4584960 \, q^{11} + 8960640 \, q^{12} + \cdots\\
\Theta_{\Lambda_{16,2,1}} &= 1 + 32 \, q^2 + 1024 \, q^3 + 8160 \, q^4 + 36864 \, q^5 + 127360 \, q^6 + 387072 \, q^7\\
 &+ 1016288 \, q^8 + 2252800 \, q^9 + 4564416 \, q^{10} + 9169920 \, q^{11} + 17395328 \, q^{12} + \cdots \\
\Theta_{\Lambda_{16,2,2}} &= 1 + 96 \, q^2 + 2048 \, q^3 + 15840 \, q^4 + 73728 \, q^5 + 259200 \, q^6 + 774144 \, q^7\\
 &+ 2003424 \, q^8 + 4505600 \, q^9 + 9269568 \, q^{10} + 18339840 \, q^{11} + 34264704 \, q^{12} +  \cdots \\
\Theta_{\Lambda_{16,2,3}} &= 1 + 224 \, q^2 + 4096 \, q^3 + 31200 \, q^4 + 147456 \, q^5 + 522880 \, q^6 + 1548288 \, q^7\\
 &+ 3977696 \, q^8 + 9011200 \, q^9 + 18679872 \, q^{10} + 36679680 \, q^{11} + 68003456 \, q^{12} +  \cdots
\end{align*}
\begin{align*}
\Theta_{{\Lambda_{16,2,1}}'} &= 1 + 32 \, q^2 + 8160 \, q^4 + 127360 \, q^6 + 1016288 \, q^8 + 4564416 \, q^{10} + 17395328 \, q^{12} + \cdots\\
\Theta_{{\Lambda_{16,2,2}}'} &= 1 + 96 \, q^2 + 15840 \, q^4 + 259200 \, q^6 + 2003424 \, q^8 + 9269568 \, q^{10} + 34264704 \, q^{12} + \cdots\\
\Theta_{{\Lambda_{16,2,3}}'} &= 1 + 224 \, q^2 + 31200 \, q^4 + 522880 \, q^6 + 3977696 \, q^8 + 18679872 \, q^{10} + 68003456 \, q^{12} + \cdots\\
\Theta_{\bm{D_{16}}^{+}} &= \Theta_{\bm{E_8} \perp \bm{E_8}}\\
 &= 1 + 480 \, q^2 + 61920 \, q^4 + 1050240 \, q^6 + 7926240 \, q^8 + 37500480 \, q^{10} + 135480960 \, q^{12} + \cdots
\end{align*}

For unimodular lattices, we have $\Theta_{\bm{D_{16}}^{+}} = \Theta_{\bm{E_8} \perp \bm{E_8}} = (\theta_3^8 - 16 \Delta_8)^2 = (\theta_3^{16} + \theta_3^{16} + \theta_3^{16})/2$ and $\Theta_{\Lambda_{16,2,3}} = \theta_3^{16} - 32 \theta_3^8 \Delta_8 = (\theta_3^{16} + \theta_3^{16} + \theta_3^{16} - 2 \theta_2^8 \theta_4^8)/2$, where $\Delta_8 = \frac{1}{16} \theta_2^4 \theta_4^4$ and $\theta_i$ for $i = 2, 3, 4$ are known as {\it Jacobi's theta functions} (See \cite{CS}). Similarly, we have $\Theta_{BW_{16}} = (\theta_3^{16} + \theta_3^{16} + \theta_3^{16} +30 \theta_2^8 \theta_3^8)/2$.\\

The following table contains the $(d, n, s, t)$-configuration of each shell of norm $m$ of the lattice:
\begin{center}
\begin{tabular}{r}
\quad\\
\hline
$m$\\
\hline
$3$\\
$4$\\
$5$\\
$6$\\
$7$\\
$8$\\
$9$\\
\hline
\end{tabular}
\quad
\begin{tabular}{cccc}
\multicolumn{4}{c}{$BW_{16}$}\\
\hline
$d$ & $n$ & $s$ & $t$\\
\hline
\ \\
$16$ & $4320$ & $6$ & $7$\\
\ \\
$16$ & $61440$ & $10$ & $7$\\
\ \\
$16$ & $522720$ & $14$ & $7$\\
\ \\
\hline
\end{tabular}
\qquad
\begin{tabular}{cccc}
\multicolumn{4}{c}{$O_{16}$}\\
\hline
$d$ & $n$ & $s$ & $t$\\
\hline
$16$ & $512$ & $4$ & $5$\\
$16$ & $4320$ & $6$ & $7$\\
$16$ & $18432$ & $8$ & $5$\\
$16$ & $61440$ & $10$ & $7$\\
$16$ & $193536$ & $12$ & $5$\\
$16$ & $522720$ & $14$ & $7$\\
$16$ & $1126400$ & $16$ & $5$\\
\hline
\end{tabular}
\end{center}\quad

\begin{center}
\begin{tabular}{r}
\quad\\
\hline
$m$\\
\hline
$2$\\
$3$\\
$4$\\
$5$\\
$6$\\
$7$\\
$8$\\
$9$\\
\hline
\end{tabular}
\quad
\begin{tabular}{cccc}
\multicolumn{4}{c}{$\Lambda_{16,2,1}$}\\
\hline
$d$ & $n$ & $s$ & $t$\\
\hline
$16$ & $32$ & $2$ & $3$\\
$16$ & $1024$ & $6$ & $5$\\
$16$ & $8160$ & $8$ & $3$\\
$16$ & $36864$ & $10$ & $5$\\
$16$ & $127360$ & $12$ & $3$\\
$16$ & $387072$ & $14$ & $5$\\
$16$ & $1016288$ & $16$ & $3$\\
$16$ & $2252800$ & $18$ & $5$\\
\hline
\end{tabular}
\quad
\begin{tabular}{cccc}
\multicolumn{4}{c}{$\Lambda_{16,2,2}$}\\
\hline
$d$ & $n$ & $s$ & $t$\\
\hline
$16$ & $96$ & $4$ & $3$\\
$16$ & $2048$ & $6$ & $5$\\
$16$ & $15840$ & $8$ & $3$\\
$16$ & $73728$ & $10$ & $5$\\
$16$ & $259200$ & $12$ & $3$\\
$16$ & $774144$ & $14$ & $5$\\
$16$ & $2003424$ & $16$ & $3$\\
\hline
\quad\\
\end{tabular}
\quad
\begin{tabular}{cccc}
\multicolumn{4}{c}{$\Lambda_{16,2,3}$}\\
\hline
$d$ & $n$ & $s$ & $t$\\
\hline
$16$ & $224$ & $4$ & $3$\\
$16$ & $4096$ & $6$ & $5$\\
$16$ & $31200$ & $8$ & $3$\\
$16$ & $147456$ & $10$ & $5$\\
$16$ & $522880$ & $12$ & $3$\\
$16$ & $1548288$ & $14$ & $5$\\
\hline
\quad\\ \quad\\
\end{tabular}
\end{center}\quad

\begin{center}
\begin{tabular}{r}
\quad\\
\hline
$m$\\
\hline
$2$\\
$4$\\
$6$\\
$8$\\
\hline
\end{tabular}
\quad
\begin{tabular}{cccc}
\multicolumn{4}{c}{${\Lambda_{16,2,1}}'$}\\
\hline
$d$ & $n$ & $s$ & $t$\\
\hline
$16$ & $32$ & $2$ & $3$\\
$16$ & $8160$ & $8$ & $3$\\
$16$ & $127360$ & $12$ & $3$\\
$16$ & $1016288$ & $16$ & $3$\\
\hline
\end{tabular}
\quad
\begin{tabular}{cccc}
\multicolumn{4}{c}{${\Lambda_{16,2,2}}'$}\\
\hline
$d$ & $n$ & $s$ & $t$\\
\hline
$16$ & $96$ & $4$ & $3$\\
$16$ & $15840$ & $8$ & $3$\\
$16$ & $259200$ & $12$ & $3$\\
$16$ & $2003424$ & $16$ & $3$\\
\hline
\end{tabular}
\quad
\begin{tabular}{cccc}
\multicolumn{4}{c}{${\Lambda_{16,2,3}}'$}\\
\hline
$d$ & $n$ & $s$ & $t$\\
\hline
$16$ & $224$ & $4$ & $3$\\
$16$ & $31200$ & $8$ & $3$\\
$16$ & $522880$ & $12$ & $3$\\
\hline
\quad\\
\end{tabular}
\end{center}\quad

\begin{center}
\begin{tabular}{r}
\quad\\
\hline
$m$\\
\hline
$2$\\
$4$\\
$6$\\
\hline
\end{tabular}
\quad
\begin{tabular}{cccc}
\multicolumn{4}{c}{$\bm{D}_{16}^{+}$}\\
\hline
$d$ & $n$ & $s$ & $t$\\
\hline
$16$ & $480$ & $4$ & $3$\\
$16$ & $61920$ & $8$ & $3$\\
$16$ & $1050240$ & $12$ & $3$\\
\hline
\end{tabular}
\quad
\begin{tabular}{cccc}
\multicolumn{4}{c}{$\bm{E}_{8} \perp \bm{E}_{8}$}\\
\hline
$d$ & $n$ & $s$ & $t$\\
\hline
$16$ & $480$ & $4$ & $3$\\
$16$ & $61920$ & $8$ & $3$\\
$16$ & $1050240$ & $12$ & $3$\\
\hline
\end{tabular}
\end{center}\quad

\newpage

\section{The laminated lattice $\Lambda_1 \simeq 2 \, \mathbb{Z}$}

We have the following obvious lattice neighboring from the laminated lattice $\Lambda_1$ to the the unimodular lattices $\mathbb{Z}$:\vspace{-0.1in}

\begin{figure}[h]
\includegraphics[width=1.45in]{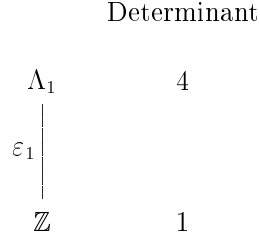}\vspace{-0.1in}
\caption{Lattice neighboring from $\Lambda_1$}\label{fig1}
\end{figure}

It is clear that $\Lambda_1^{\sharp} = \Lambda_1 \cup (\varepsilon_1 + \Lambda_1) \cup (\frac{1}{2} \varepsilon_1 + \Lambda_1) \cup (- \frac{1}{2} \varepsilon_1 + \Lambda_1)$, and unique lattice which include $\Lambda_1$ is $\Lambda_1 \cup (\varepsilon_1 + \Lambda_1) = \mathbb{Z}$.\\

\subsection{Properties of lattices}

The theta series of each lattice have the following form:
\begin{align*}
\Theta_{\Lambda_1} &= 1 + 2 \, q^{4} + 2 \, q^{16} + 2 \, q^{36} + 2 \, q^{64} + 2 \, q^{100} + \cdots\\
\Theta_{\mathbb{Z}} &= 1 + 2 \, q + 2 \, q^4 + 2 \, q^9 + 2 \, q^{16} + 2 \, q^{25} + \cdots
\end{align*}

On spherical designs, it is obvious that every nonempty shell of each lattice satisfies the equation (\ref{eq-design}) for any polynomials of all degree $t$.\\

\section{The laminated lattice $\Lambda_2 \simeq \sqrt{2} \, \bm{A}_2$}

We have the following lattice neighboring from the laminated lattice $\Lambda_2$:\vspace{-0.1in}

\begin{figure}[h]
\includegraphics[width=1.55in]{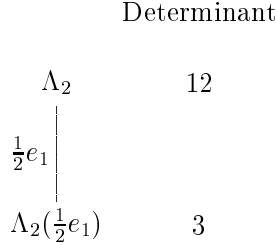}\vspace{-0.1in}
\caption{Lattice neighboring from $\Lambda_2$}\label{fig2}
\end{figure}

Here, a Gram matrix of $\Lambda_2$ is ${}^t \! (e_1, e_2) \cdot (e_1, e_2) = 2 \left[ \begin{array}{cc} 2 & -1\\ -1 & 2 \end{array}\right]$, and we have $\Lambda_2 (\frac{1}{2} e_1) \simeq \mathbb{Z} \perp O_1$.

We can show the following fact:

\begin{theorem}
The integral lattices which include $\Lambda_2$ are isometric to one of the above two lattices $($in Figure $\ref{fig2}$.$)$, namely $\Lambda_2$ and $\Lambda_2 (\frac{1}{2} e_1)$.
\end{theorem}

It is easy to classify the above lattices. We have the following:

\begin{proposition}
The $12$ classes of $\Lambda_2^{\sharp} / \Lambda_2$ consisit of the following classes:\\

\begin{center}
\begin{tabular}{ccc}
\hline
\# of class & norm of c.l. & order\\
\hline
$1$ & $0$ & $1$\\
$3 \times 2$ & $1/3$ & $6$\\
$3$ & $1$ & $2$\\
$1 \times 2$ & $4/3$ & $3$\\
\hline
\end{tabular}
\end{center}\quad

\noindent
where $\times 2$ means that classes are represented by some pairs of vectors as $v$ and $- v$, and `c.l.' means coset leader.
\end{proposition}

Thus, to classify such lattices, we have only to consider the neighboring with the each vector of norm $1$. We can easily check by computer search. In fact, we have $4$ integral lattices, one of which is $\Lambda_2$ and the other three of which are isometric to $\Lambda_2 (\frac{1}{2} e_1)$.\\

\subsection{Properties of lattices}

The theta series of each lattice have the following form:
\begin{align*}
\Theta_{\Lambda_2} &= 1 + 6 \, q^4 + 6 \, q^{12} + \cdots\\
\Theta_{\Lambda_2 (\frac{1}{2} e_1)} &= 1 + 2 \, q + 2 \, q^3 + 6 \, q^4 + 4 \, q^7 + 2 \, q^9 + 6 \, q^{12} + \cdots
\end{align*}

The following table contains the $(d, n, s, t)$-configuration of each shell of norm $m$ of the lattice:
\begin{center}
\begin{tabular}{r}
\quad\\
\hline
$m$\\
\hline
$1$\\
$2$\\
$3$\\
$4$\\
$5$\\
$6$\\
$7$\\
$8$\\
$9$\\
$10$\\
$11$\\
$12$\\
\hline
\end{tabular}
\quad
\begin{tabular}{cccc}
\multicolumn{4}{c}{$\Lambda_2$}\\
\hline
$d$ & $n$ & $s$ & $t$\\
\hline
\ \\ \ \\ \ \\
$2$ & $6$ & $3$ & $5$\\
\ \\ \ \\ \ \\ \ \\ \ \\ \ \\ \ \\
$2$ & $6$ & $3$ & $5$\\
\hline
\end{tabular}
\qquad
\begin{tabular}{cccc}
\multicolumn{4}{c}{$\Lambda_2 (\frac{1}{2} e_1)$}\\
\hline
$d$ & $n$ & $s$ & $t$\\
\hline
$1$ & $2$ & $1$ & $-$\\
\ \\
$1$ & $2$ & $1$ & $-$\\
$2$ & $6$ & $3$ & $5$\\
\ \\ \ \\
$2$ & $4$ & $3$ & $1$\\
\ \\
$1$ & $2$ & $1$ & $-$\\
\ \\ \ \\
$2$ & $6$ & $3$ & $5$\\
\hline
\end{tabular}
\end{center}\quad\\

\section{The laminated lattice $\Lambda_3 \simeq \sqrt{2} \, \bm{A}_3$}

We have the following lattice neighboring from the laminated lattice $\Lambda_3$:\vspace{-0.1in}

\begin{figure}[h]
\includegraphics[width=3in]{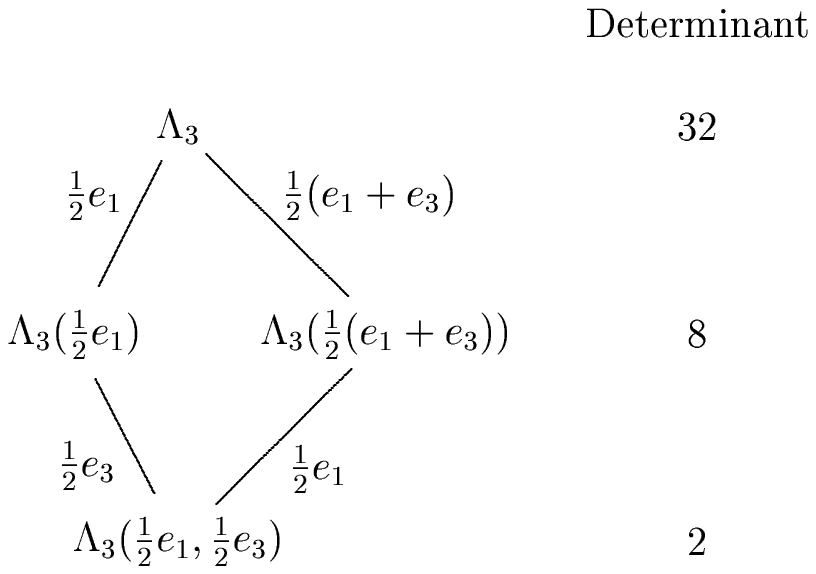}\vspace{-0.1in}
\caption{Lattice neighboring from $\Lambda_3$}\label{fig3}
\end{figure}

Here, a Gram matrix of $\Lambda_3$ is
\begin{equation*}
{}^t \! (e_1, e_2, e_3) \cdot (e_1, e_2, e_3) = 2 {\small \left[ \begin{array}{ccc} 2 & -1 & 0\\ -1 & 2 & -1\\ 0 & -1 & 2 \end{array} \right]},
\end{equation*}
 and we have $\Lambda_3 (\frac{1}{2} (e_1 + e_3)) \simeq (\bm{A}_1)^3$ and $\Lambda_3 (\frac{1}{2} e_1, \frac{1}{2} e_3) \simeq \mathbb{Z}^2 \perp \bm{A}_1$.

We can show the following fact:
\begin{theorem}
The integral lattices which include $\Lambda_3$ are isometric to one of the above four lattices $($in Figure $\ref{fig3}$.$)$, namely $\Lambda_3$, $\Lambda_3 (\frac{1}{2} e_1)$, $\Lambda_3 (\frac{1}{2} (e_1 + e_3))$, and $\Lambda_3 (\frac{1}{2} e_1, \frac{1}{2} e_3)$.
\end{theorem}

\newpage

It is easy to classify the above lattices. We have the following:
\begin{proposition}
The $32$ classes of $\Lambda_3^{\sharp} / \Lambda_3$ consisit of the following classes:\\

\begin{center}
\begin{tabular}{ccc}
\hline
\# of class & norm of c.l. & order\\
\hline
$1$ & $0$ & $1$\\
$4 \times 2$ & $3/8$ & $8$\\
$3 \times 2$ & $1/2$ & $4$\\
$6$ & $1$ & $2$\\
$4 \times 2$ & $11/8$ & $8$\\
$1 \times 2$ & $3/2$ & $4$\\
$1$ & $2$ & $2$\\
\hline
\end{tabular}
\end{center}\quad
\end{proposition}

Thus, to classify such lattices, we have only to consider the neighboring with the each vector of norm $1$ and $2$. We can easily check by computer search. In fact, we have $12$ integral lattices from the combinations of the above classes. The following table contains the number of integral lattices which are isometric to each lattice:

\begin{center}
\begin{tabular}{c|ccc|ccc|c}
\cline{1-2} \cline{4-5} \cline{7-8}
$\Lambda_3$ & $1$ & \qquad & $\Lambda_3 (\frac{1}{2} e_1)$ & $6$ & \qquad & $\Lambda_3 (\frac{1}{2} e_1, \frac{1}{2} e_3)$ & $4$\\
  &&& $\Lambda_3 (\frac{1}{2} (e_1 + e_3))$ & $1$ &&&\\
\cline{1-2} \cline{4-5} \cline{7-8}
\end{tabular}
\end{center}\quad

\subsection{Properties of lattices}

The theta series of each lattice have the following form:
\begin{align*}
\Theta_{\Lambda_3} &= 1 + 12 \, q^4 + 6 \, q^8 + 24 \, q^{12} + \cdots\\
\Theta_{\Lambda_3 (\frac{1}{2} e_1)} &= 1 + 2 \, q + 4 \, q^3 + 12 \, q^4 + 4 \, q^5 + 8 \, q^7\\
 &+ 6 \, q^8 + 6 \, q^9 + 4 \, q^{11} + 24 \, q^{12} + \cdots\\
\Theta_{\Lambda_3 (\frac{1}{2} (e_1 + e_3))} &= 1 + 6 \, q^2 + 12 \, q^4 + 8 \, q^6 + 6 \, q^8 + 24 \, q^{10} + 24 \, q^{12} + \cdots\\
\Theta_{\Lambda_3 (\frac{1}{2} e_1, \frac{1}{2} e_3)} &= 1 + 4 \, q + 6 \, q^2 + 8 \, q^3 + 12 \, q^4 + 8 \, q^5 + 8 \, q^6 + 16 \, q^7\\
 &+ 6 \, q^8 + 12 \, q^9 + 24 \, q^{10} + 8 \, q^{11} + 24 \, q^{12} + \cdots
\end{align*}

The following table contains the $(d, n, s, t)$-configuration of each shell of norm $m$ of the lattice:
\begin{center}
\begin{tabular}{r}
\quad\\
\hline
$m$\\
\hline
$1$\\
$2$\\
$3$\\
$4$\\
$5$\\
$6$\\
$7$\\
$8$\\
$9$\\
$10$\\
$11$\\
$12$\\
\hline
\end{tabular}
\quad
\begin{tabular}{cccc}
\multicolumn{4}{c}{$\Lambda_3$}\\
\hline
$d$ & $n$ & $s$ & $t$\\
\hline
\ \\ \ \\ \ \\
$3$ & $12$ & $4$ & $3$\\
\ \\ \ \\ \ \\
$3$ & $6$ & $2$ & $3$\\
\ \\ \ \\ \ \\
$3$ & $24$ & $11$ & $3$\\
\hline
\end{tabular}
\qquad
\begin{tabular}{cccc}
\multicolumn{4}{c}{$\Lambda_3 (\frac{1}{2} e_1)$}\\
\hline
$d$ & $n$ & $s$ & $t$\\
\hline
$1$ & $2$ & $1$ & $-$\\
\ \\
$2$ & $4$ & $3$ & $1$\\
$3$ & $12$ & $4$ & $3$\\
$2$ & $4$ & $3$ & $1$\\
\ \\
$3$ & $8$ & $7$ & $1$\\
$3$ & $6$ & $2$ & $3$\\
$2$ & $6$ & $5$ & $1$\\
\ \\
$2$ & $4$ & $3$ & $1$\\
$3$ & $24$ & $11$ & $3$\\
\hline
\end{tabular}
\qquad
\begin{tabular}{cccc}
\multicolumn{4}{c}{$\Lambda_3 (\frac{1}{2} (e_1 + e_3))$}\\
\hline
$d$ & $n$ & $s$ & $t$\\
\hline
\ \\
$3$ & $6$ & $2$ & $3$\\
\ \\
$3$ & $12$ & $4$ & $3$\\
\ \\
$3$ & $8$ & $3$ & $3$\\
\ \\
$3$ & $6$ & $2$ & $3$\\
\ \\
$3$ & $24$ & $10$ & $3$\\
\ \\
$3$ & $24$ & $11$ & $3$\\
\hline
\end{tabular}
\qquad
\begin{tabular}{cccc}
\multicolumn{4}{c}{$\Lambda_3 (\frac{1}{2} e_1, \frac{1}{2} e_3)$}\\
\hline
$d$ & $n$ & $s$ & $t$\\
\hline
$2$ & $4$ & $2$ & $3$\\
$3$ & $6$ & $2$ & $3$\\
$3$ & $8$ & $5$ & $1$\\
$3$ & $12$ & $4$ & $3$\\
$2$ & $8$ & $6$ & $3$\\
$3$ & $8$ & $3$ & $3$\\
$3$ & $16$ & $11$ & $1$\\
$3$ & $6$ & $2$ & $3$\\
$3$ & $12$ & $8$ & $1$\\
$3$ & $24$ & $10$ & $3$\\
$3$ & $8$ & $5$ & $1$\\
$3$ & $24$ & $11$ & $3$\\
\hline
\end{tabular}
\end{center}

\newpage

\section{The laminated lattice $\Lambda_4 \simeq \sqrt{2} \, \bm{D}_4$}

We have the following lattice neighboring from the laminated lattice $\Lambda_4$:\vspace{-0.1in}

\begin{figure}[h]
\includegraphics[width=2.4in]{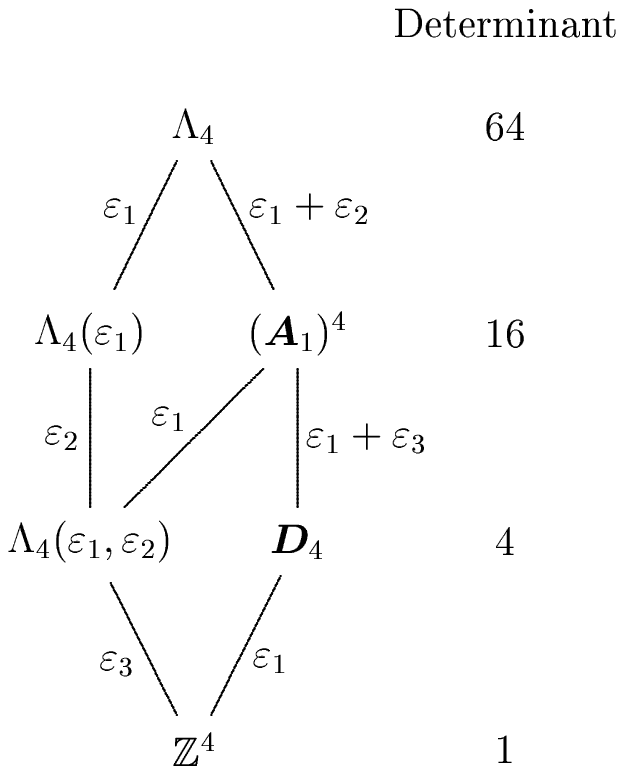}\vspace{-0.1in}
\caption{Lattice neighboring from $\Lambda_4$}\label{fig4}
\end{figure}

Here, for the above neighboring, we have the following generator matrix of $\Lambda_4$:
{\small \begin{equation*}
\left(\begin{array}{cccc}
2 & 0 & 0 & 0\\
0 & 2 & 0 & 0\\
0 & 0 & 2 & 0\\
1 & 1 & 1 & 1\\
\end{array}\right)
\end{equation*}}
We also have $\Lambda_4(\varepsilon_1, \varepsilon_2) \simeq \mathbb{Z}^2 \perp (\bm{A}_1)^2$.

We can show the following fact:
\begin{theorem}
The integral lattices which include $\Lambda_4$ are isometric to one of the above six lattices $($in Figure $\ref{fig4}$.$)$, namely $\Lambda_4$, $\Lambda_4 (\varepsilon_1)$, $\Lambda_4 (\varepsilon_1, \varepsilon_2)$, $(\bm{A}_1)^{4}$, $\bm{D}_4$ and $\mathbb{Z}^4$.
\end{theorem}

It is easy to classify the above lattices. We have the following:
\begin{proposition}
The $64$ classes of $\Lambda_4^{\sharp} / \Lambda_4$ consisit of the following classes:
\begin{center}
\begin{tabular}{ccc}
\hline
\# of class & norm of c.l. & order\\
\hline
$1$ & $0$ & $1$\\
$12 \times 2$ & $1/2$ & $4$\\
$12$ & $1$ & $2$\\
$12 \times 2$ & $3/2$ & $4$\\
$3$ & $2$ & $2$\\
\hline
\end{tabular}
\end{center}
\end{proposition}

Thus, to classify such lattices, we have only to consider the neighboring with the each vector of norm $1$ and $2$. We can easily check by computer search. In fact, we have $38$ integral lattices from the combinations of the above classes. The following table contains the number of integral lattices which are isometric to each lattice:

\begin{center}
\begin{tabular}{c|ccc|ccc|ccc|c}
\cline{1-2} \cline{4-5} \cline{7-8} \cline{10-11}
$\Lambda_4$ & $1$ & \qquad & $\Lambda_4 (\varepsilon_1)$ & $12$ & \qquad & $\Lambda_4 (\varepsilon_1, \varepsilon_2)$ & $18$ & \qquad & $\mathbb{Z}^4$ & $3$\\
  &&& $(\bm{A}_1)^{4}$ & $3$ && $\bm{D}_4$ & $1$ &&&\\
\cline{1-2} \cline{4-5} \cline{7-8} \cline{10-11}
\end{tabular}
\end{center}\quad

\subsection{Properties of lattices}

The theta series of each lattice have the following form:
\begin{align*}
\Theta_{\Lambda_4} &= 1 + 24 \, q^4 + 24 \, q^8 + 96 \, q^{12} + \cdots\\
\Theta_{(\bm{A}_1)^{4}} &= 1 + 8 \, q^2 + 24 \, q^4 + 32 \, q^6 + 24 \, q^8 + 48 \, q^{10} + 96 \, q^{12} + \cdots\\
\Theta_{\bm{D}_4} &= 1 + 24 \, q^2 + 24 \, q^4 + 96 \, q^6 + 24 \, q^8 + 144 \, q^{10} + 96 \, q^{12} + \cdots\\
\quad\\
\Theta_{\Lambda_4 (\varepsilon_1)} &= 1 + 2 \, q + 8 \, q^3 + 24 \, q^4 + 12 \, q^5 + 16 \, q^7\\
 &+ 24 \, q^8 + 26 \, q^9 + 24 \, q^{11} + 96 \, q^{12} + \cdots\\
\Theta_{\Lambda_4 \varepsilon_1, \varepsilon_2)} &= 1 + 4 \, q + 8 \, q^2 + 16 \, q^3 + 24 \, q^4 + 24 \, q^5 + 32 \, q^6 + 32 \, q^7\\
 &+ 24 \, q^8 + 52 \, q^9 + 48 \, q^{10} + 48 \, q^{11} + 96 \, q^{12} + \cdots\\
\Theta_{\mathbb{Z}^4} &= 1 + 8 \, q + 24 \, q^2 + 32 \, q^3 + 24 \, q^4 + 48 \, q^5 + 96 \, q^6 + 64 \, q^7\\
 &+ 24 \, q^8 + 104 \, q^9 + 144 \, q^{10} + 96 \, q^{11} + 96 \, q^{12} + \cdots
\end{align*}

The following table contains the $(d, n, s, t)$-configuration of each shell of norm $m$ of the lattice:
\begin{center}
\begin{tabular}{r}
\quad\\
\hline
$m$\\
\hline
$2$\\
$4$\\
$6$\\
$8$\\
$10$\\
$12$\\
\hline
\end{tabular}
\quad
\begin{tabular}{cccc}
\multicolumn{4}{c}{$\Lambda_4$}\\
\hline
$d$ & $n$ & $s$ & $t$\\
\hline
\ \\
$4$ & $24$ & $4$ & $5$\\
\ \\
$4$ & $24$ & $4$ & $5$\\
\ \\
$4$ & $96$ & $12$ & $5$\\
\hline
\end{tabular}
\qquad
\begin{tabular}{cccc}
\multicolumn{4}{c}{$(\bm{A}_1)^{4}$}\\
\hline
$d$ & $n$ & $s$ & $t$\\
\hline
$4$ & $8$ & $2$ & $3$\\
$4$ & $24$ & $4$ & $5$\\
$4$ & $32$ & $6$ & $3$\\
$4$ & $24$ & $4$ & $5$\\
$4$ & $48$ & $10$ & $3$\\
$4$ & $96$ & $12$ & $5$\\
\hline
\end{tabular}
\qquad
\begin{tabular}{cccc}
\multicolumn{4}{c}{$\bm{D}_4$}\\
\hline
$d$ & $n$ & $s$ & $t$\\
\hline
$4$ & $24$ & $4$ & $5$\\
$4$ & $24$ & $4$ & $5$\\
$4$ & $96$ & $12$ & $5$\\
$4$ & $24$ & $4$ & $5$\\
$4$ & $144$ & $20$ & $5$\\
$4$ & $96$ & $12$ & $5$\\
\hline
\end{tabular}
\end{center}\quad

\begin{center}
\begin{tabular}{r}
\quad\\
\hline
$m$\\
\hline
$1$\\
$2$\\
$3$\\
$4$\\
$5$\\
$6$\\
$7$\\
$8$\\
$9$\\
$10$\\
$11$\\
$12$\\
\hline
\end{tabular}
\qquad
\begin{tabular}{cccc}
\multicolumn{4}{c}{$\Lambda_4 (\varepsilon_1)$}\\
\hline
$d$ & $n$ & $s$ & $t$\\
\hline
$1$ & $2$ & $1$ & $-$\\
\ \\
$3$ & $8$ & $3$ & $3$\\
$4$ & $24$ & $4$ & $5$\\
$4$ & $12$ & $5$ & $1$\\
\ \\
$4$ & $16$ & $7$ & $1$\\
$4$ & $24$ & $4$ & $5$\\
$4$ & $26$ & $9$ & $1$\\
\ \\
$3$ & $24$ & $9$ & $3$\\
$4$ & $96$ & $12$ & $5$\\
\hline
\end{tabular}
\qquad
\begin{tabular}{cccc}
\multicolumn{4}{c}{$\Lambda_4 (\varepsilon_1, \varepsilon_2)$}\\
\hline
$d$ & $n$ & $s$ & $t$\\
\hline
$2$ & $4$ & $2$ & $3$\\
$4$ & $8$ & $2$ & $3$\\
$4$ & $16$ & $6$ & $1$\\
$4$ & $24$ & $4$ & $5$\\
$4$ & $24$ & $10$ & $1$\\
$4$ & $32$ & $6$ & $3$\\
$4$ & $32$ & $14$ & $1$\\
$4$ & $24$ & $4$ & $5$\\
$4$ & $52$ & $18$ & $1$\\
$4$ & $48$ & $10$ & $3$\\
$4$ & $48$ & $20$ & $1$\\
$4$ & $96$ & $12$ & $5$\\
\hline
\end{tabular}
\qquad
\begin{tabular}{cccc}
\multicolumn{4}{c}{$\mathbb{Z}^4$}\\
\hline
$d$ & $n$ & $s$ & $t$\\
\hline
$4$ & $8$ & $2$ & $3$\\
$4$ & $24$ & $4$ & $5$\\
$4$ & $32$ & $6$ & $3$\\
$4$ & $24$ & $4$ & $5$\\
$4$ & $48$ & $10$ & $3$\\
$4$ & $96$ & $12$ & $5$\\
$4$ & $64$ & $14$ & $3$\\
$4$ & $24$ & $4$ & $5$\\
$4$ & $104$ & $18$ & $3$\\
$4$ & $144$ & $20$ & $5$\\
$4$ & $96$ & $20$ & $3$\\
$4$ & $96$ & $12$ & $5$\\
\hline
\end{tabular}
\end{center}

\newpage

\section{The laminated lattice $\Lambda_5 \simeq \sqrt{2} \, \bm{D}_5$}

We have the following lattice neighboring from the laminated lattice $\Lambda_5$:\vspace{-0.1in}

\begin{figure}[h]
\includegraphics[width=5.3in]{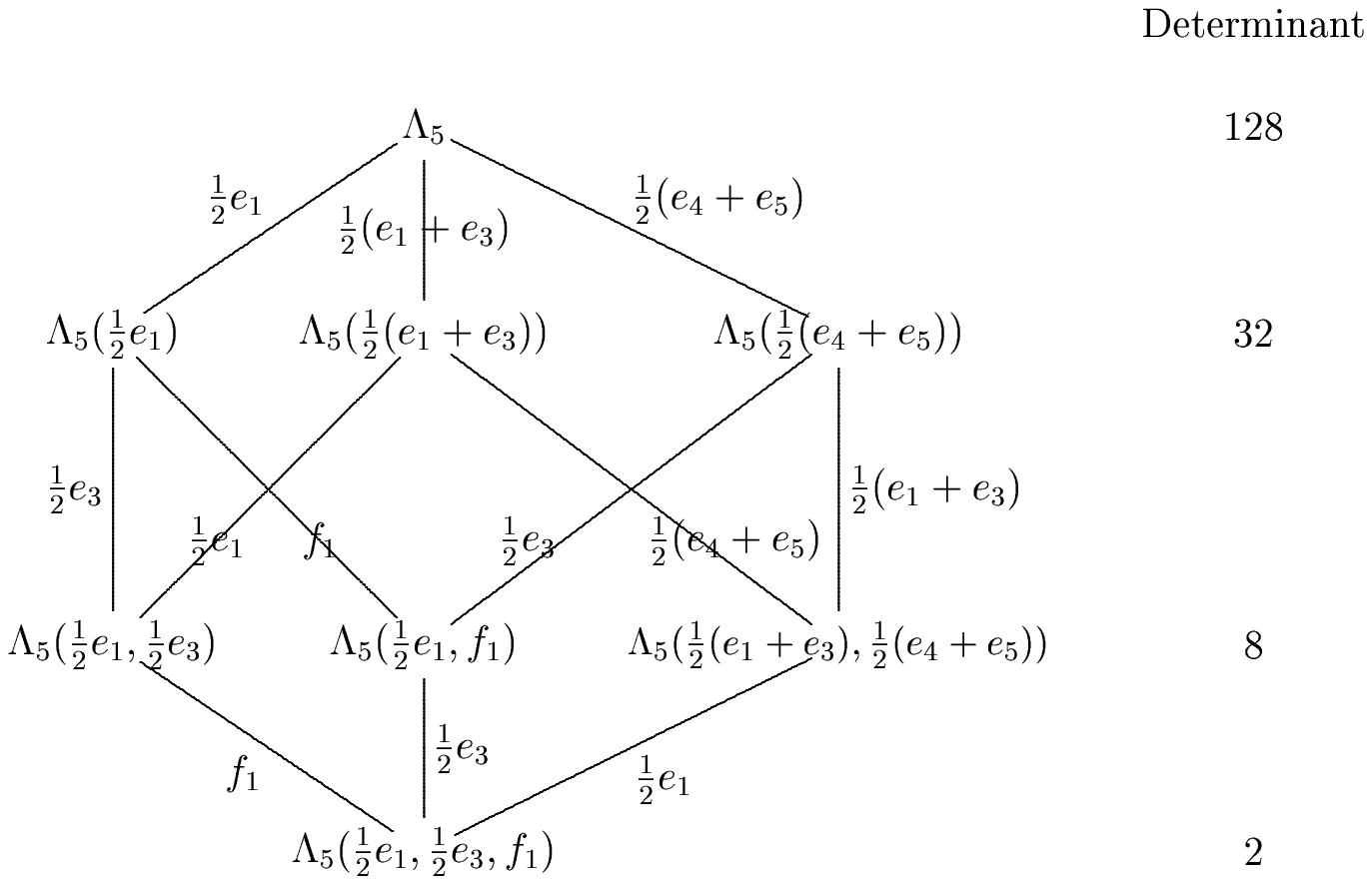}\vspace{-0.1in}
\caption{Lattice neighboring from $\Lambda_5$}\label{fig5}
\end{figure}

Here, a Gram matrix of $\Lambda_5$ is
\begin{equation*}
{}^t \! (e_1, \ldots, e_5) \cdot (e_1, \ldots, e_5) = 2
 {\small \left[ \begin{array}{ccccc} 2 & -1 & 0 & 0 & 0\\ -1 & 2 & -1 & 0 & 0\\ 0 & -1 & 2 & -1 & -1\\ 0 & 0 & -1 & 2 & 0\\ 0 & 0 & -1 & 0 & 2 \end{array} \right]},
\end{equation*}
we denote $f_1 := \frac{1}{2} (e_1 + 2 e_2 + 2 e_3 + e_4 + e_5)$, and we have $\Lambda_5 (\frac{1}{2} (e_4 + e_5)) \simeq (\bm{A}_1)^5$, $\Lambda_5 (\frac{1}{2} e_1, f_1) \simeq \mathbb{Z}^2 \perp (\bm{A}_1)^3$, $\Lambda_5 (\frac{1}{2} (e_1 + e_3), \frac{1}{2} (e_4 + e_5)) \simeq \bm{A}_1 \perp \bm{D}_4$, and $\Lambda_5 (\frac{1}{2} e_1, \frac{1}{2} e_3, f_1) \simeq \mathbb{Z}^4 \perp \bm{A}_1$.

We can show the following fact:
\begin{theorem}
The integral lattices which include $\Lambda_5$ are isometric to one of the above eight lattices $($in Figure $\ref{fig5}$.$)$, namely $\Lambda_5$, $\Lambda_5 (\frac{1}{2} e_1)$, $\Lambda_5 (\frac{1}{2} (e_1 + e_3))$, $\Lambda_5 (\frac{1}{2} (e_4 + e_5))$, $\Lambda_5 (\frac{1}{2} e_1, \frac{1}{2} e_3)$, $\Lambda_5 (\frac{1}{2} e_1, f_1)$, $\Lambda_5 (\frac{1}{2} (e_1 + e_3), \frac{1}{2} (e_4 + e_5))$, and $\Lambda_5 (\frac{1}{2} e_1, \frac{1}{2} e_3, f_1)$.
\end{theorem}

It is easy to classify the above lattices. We have the following:
\begin{proposition}
The $128$ classes of $\Lambda_5^{\sharp} / \Lambda_5$ consisit of the following classes:\\

\begin{center}
\begin{tabular}{ccc}
\hline
\# of class & norm of c.l. & order\\
\hline
$1$ & $0$ & $1$\\
$5 \times 2$ & $1/2$ & $4$\\
$16 \times 2$ & $5/8$ & $8$\\
$20$ & $1$ & $2$\\
$10 \times 2$ & $3/2$ & $4$\\
$16 \times 2$ & $13/8$ & $8$\\
$11$ & $2$ & $2$\\
$1 \times 2$ & $5/2$ & $4$\\
\hline
\end{tabular}
\end{center}\quad
\end{proposition}

Thus, to classify such lattices, we have only to consider the neighboring with the each vector of norm $1$ and $2$. We can easily check by computer search. In fact, we have $122$ integral lattices from the combinations of the above classes. The following table contains the number of integral lattices which are isometric to each lattice:\\

\begin{center}
\begin{tabular}{c|ccc|ccc|ccc|c}
\cline{1-2} \cline{4-5} \cline{7-8} \cline{10-11}
$\Lambda_5$ & $1$ & \qquad & $\Lambda_5 (\frac{1}{2} e_1)$ & $20$ & \qquad & $\Lambda_5 (\frac{1}{2} e_1, \frac{1}{2} e_3)$ & $60$ & \qquad & $\Lambda_5 (\frac{1}{2} e_1, \frac{1}{2} e_3, f_1)$ & $15$\\
  &&& $\Lambda_5 (\frac{1}{2} (e_1 + e_3))$ & $10$ && $\Lambda_5 (\frac{1}{2} e_1, f_1)$ & $10$ &&&\\
  &&& $\Lambda_5 (\frac{1}{2} (e_4 + e_5))$ & $1$ && {\small $\Lambda_5 (\frac{1}{2} (e_1 + e_3), \frac{1}{2} (e_4 + e_5))$} & $5$ &&&\\
\cline{1-2} \cline{4-5} \cline{7-8} \cline{10-11}
\end{tabular}
\end{center}\quad

\subsection{Properties of lattices}

The theta series of each lattice have the following form:
\begin{align*}
\Theta_{\Lambda_5} &= 1 + 40 \, q^4 + 90 \, q^8 + 240 \, q^{12} + \cdots\\
\Theta_{\Lambda_5 (\frac{1}{2} (e_1 + e_3))} &= 1 + 8 \, q^2 + 56 \, q^4 + 32 \, q^6 + 154 \, q^8 + 64 \, q^{10} + 336 \, q^{12} + \cdots\\
\Theta_{\Lambda_5 (\frac{1}{2} (e_4 + e_5))} &= 1 + 10 \, q^2 + 40 \, q^4 + 80 \, q^6 + 90 \, q^8 + 112 \, q^{10} + 240 \, q^{12} + \cdots\\
\Theta_{\Lambda_5 (\frac{1}{2} (e_1 + e_3), \frac{1}{2} (e_4 + e_5))} &= 1 + 26 \, q^2 + 72 \, q^4 + 144 \, q^6 + 218 \, q^8 + 240 \, q^{10} + 432 \, q^{12} + \cdots\\
\quad\\
\Theta_{\Lambda_5 (\frac{1}{2} e_1)} &= 1 + 2 \, q + 12 \, q^3 + 40 \, q^4 + 28 \, q^5 + 40 \, q^7\\
 &+ 90 \, q^8 + 62 \, q^9 + 92 \, q^{11} + 240 \, q^{12} + \cdots\\
\Theta_{\Lambda_5 (\frac{1}{2} e_1, \frac{1}{2} e_3)} &= 1 + 4 \, q + 8 \, q^2 + 24 \, q^3 + 56 \, q^4 + 56 \, q^5 + 32 \, q^6 + 80 \, q^7\\
 &+ 154 \, q^8 + 124 \, q^9 + 64 \, q^{10} + 184 \, q^{11} + 336 \, q^{12} + \cdots\\
\Theta_{\Lambda_5 (\frac{1}{2} e_1, f_1)} &= 1 + 4 \, q + 10 \, q^2 + 24 \, q^3 + 40 \, q^4 + 56 \, q^5 + 80 \, q^6 + 80 \, q^7\\
 &+ 90 \, q^8 + 124 \, q^9 + 112 \, q^{10} + 184 \, q^{11} + 240 \, q^{12} + \cdots\\
\Theta_{\Lambda_5 (\frac{1}{2} e_1, \frac{1}{2} e_3, f_1)} &= 1 + 8 \, q + 26 \, q^2 + 48 \, q^3 + 72 \, q^4 + 112 \, q^5 + 144 \, q^6 + 160 \, q^7\\
 &+ 218 \, q^8 + 248 \, q^9 + 240 \, q^{10} + 368 \, q^{11} + 432 \, q^{12} + \cdots
\end{align*}\quad

The following table contains the $(d, n, s, t)$-configuration of each shell of norm $m$ of the lattice:
\begin{center}
\begin{tabular}{r}
\quad\\
\hline
$m$\\
\hline
$2$\\
$4$\\
$6$\\
$8$\\
$10$\\
$12$\\
\hline
\end{tabular}
\quad
\begin{tabular}{cccc}
\multicolumn{4}{c}{$\Lambda_5$}\\
\hline
$d$ & $n$ & $s$ & $t$\\
\hline
\ \\
$5$ & $40$ & $4$ & $3$\\
\ \\
$5$ & $90$ & $8$ & $3$\\
\ \\
$5$ & $240$ & $12$ & $3$\\
\hline
\end{tabular}
\qquad
\begin{tabular}{cccc}
\multicolumn{4}{c}{$\Lambda_5 (\frac{1}{2} (e_1 + e_3))$}\\
\hline
$d$ & $n$ & $s$ & $t$\\
\hline
$4$ & $8$ & $2$ & $3$\\
$5$ & $56$ & $8$ & $1$\\
$4$ & $32$ & $6$ & $3$\\
$5$ & $154$ & $16$ & $1$\\
$5$ & $64$ & $10$ & $3$\\
$5$ & $336$ & $24$ & $1$\\
\hline
\end{tabular}
\qquad
\begin{tabular}{cccc}
\multicolumn{4}{c}{$\Lambda_5 (\frac{1}{2} (e_4 + e_5))$}\\
\hline
$d$ & $n$ & $s$ & $t$\\
\hline
$5$ & $10$ & $2$ & $3$\\
$5$ & $40$ & $4$ & $3$\\
$5$ & $80$ & $6$ & $3$\\
$5$ & $90$ & $8$ & $3$\\
$5$ & $112$ & $10$ & $3$\\
$5$ & $240$ & $12$ & $3$\\
\hline
\end{tabular}
\qquad
\begin{tabular}{cccc}
\multicolumn{4}{c}{\tiny $\Lambda_5 (\frac{1}{2} (e_1 + e_3), \frac{1}{2} (e_4 + e_5))$}\\
\hline
\, $d$ \, & \, $n$ \, & \, $s$ \, & \, $t$ \, \\
\hline
$5$ & $26$ & $4$ & $1$\\
$5$ & $72$ & $8$ & $1$\\
$5$ & $144$ & $12$ & $1$\\
$5$ & $218$ & $16$ & $1$\\
$5$ & $240$ & $20$ & $3$\\
$5$ & $432$ & $24$ & $1$\\
\hline
\end{tabular}
\end{center}\quad

\begin{center}
\begin{tabular}{r}
\quad\\
\hline
$m$\\
\hline
$1$\\
$2$\\
$3$\\
$4$\\
$5$\\
$6$\\
$7$\\
$8$\\
$9$\\
$10$\\
$11$\\
$12$\\
\hline
\end{tabular}
\qquad
\begin{tabular}{cccc}
\multicolumn{4}{c}{$\Lambda_5 (\frac{1}{2} e_1)$}\\
\hline
$d$ & $n$ & $s$ & $t$\\
\hline
$1$ & $2$ & $1$ & $-$\\
\ \\
$4$ & $12$ & $3$ & $1$\\
$5$ & $40$ & $4$ & $3$\\
$5$ & $28$ & $5$ & $1$\\
\ \\
$5$ & $40$ & $7$ & $1$\\
$5$ & $90$ & $8$ & $3$\\
$5$ & $62$ & $9$ & $1$\\
\ \\
$5$ & $92$ & $11$ & $1$\\
$5$ & $240$ & $12$ & $3$\\
\hline
\end{tabular}
\qquad
\begin{tabular}{cccc}
\multicolumn{4}{c}{$\Lambda_5 (\frac{1}{2} e_1, \frac{1}{2} e_3)$}\\
\hline
$d$ & $n$ & $s$ & $t$\\
\hline
$2$ & $4$ & $2$ & $3$\\
$4$ & $8$ & $2$ & $3$\\
$5$ & $24$ & $6$ & $1$\\
$5$ & $56$ & $8$ & $1$\\
$5$ & $24$ & $10$ & $1$\\
$4$ & $32$ & $6$ & $3$\\
$5$ & $80$ & $14$ & $1$\\
$5$ & $154$ & $16$ & $1$\\
$5$ & $124$ & $18$ & $1$\\
$5$ & $64$ & $10$ & $3$\\
$5$ & $184$ & $22$ & $1$\\
$5$ & $336$ & $24$ & $1$\\
\hline
\end{tabular}
\qquad
\begin{tabular}{cccc}
\multicolumn{4}{c}{$\Lambda_5 (\frac{1}{2} e_1, f_1)$}\\
\hline
$d$ & $n$ & $s$ & $t$\\
\hline
$2$ & $4$ & $2$ & $3$\\
$5$ & $10$ & $2$ & $3$\\
$5$ & $24$ & $6$ & $1$\\
$5$ & $40$ & $4$ & $3$\\
$5$ & $56$ & $10$ & $1$\\
$5$ & $80$ & $6$ & $3$\\
$5$ & $80$ & $14$ & $1$\\
$5$ & $90$ & $8$ & $3$\\
$5$ & $124$ & $18$ & $1$\\
$5$ & $112$ & $10$ & $3$\\
$5$ & $184$ & $22$ & $1$\\
$5$ & $240$ & $12$ & $3$\\
\hline
\end{tabular}
\qquad
\begin{tabular}{cccc}
\multicolumn{4}{c}{\small $\Lambda_5 (\frac{1}{2} e_1, \frac{1}{2} e_3, f_1)$}\\
\hline
$d$ & $n$ & $s$ & $t$\\
\hline
$4$ & $8$ & $2$ & $3$\\
$5$ & $26$ & $4$ & $1$\\
$5$ & $48$ & $6$ & $1$\\
$5$ & $72$ & $8$ & $1$\\
$5$ & $112$ & $10$ & $1$\\
$5$ & $144$ & $12$ & $1$\\
$5$ & $160$ & $14$ & $1$\\
$5$ & $218$ & $16$ & $1$\\
$5$ & $248$ & $18$ & $1$\\
$5$ & $240$ & $20$ & $3$\\
$5$ & $368$ & $22$ & $1$\\
$5$ & $432$ & $24$ & $1$\\
\hline
\end{tabular}
\end{center}

\newpage

\section{The laminated lattice $\Lambda_6 \simeq \sqrt{2} \, \bm{E}_6$}

We have the following lattice neighboring from the laminated lattice $\Lambda_6$:\vspace{-0.1in}

\begin{figure}[h]
\includegraphics[width=3.4in]{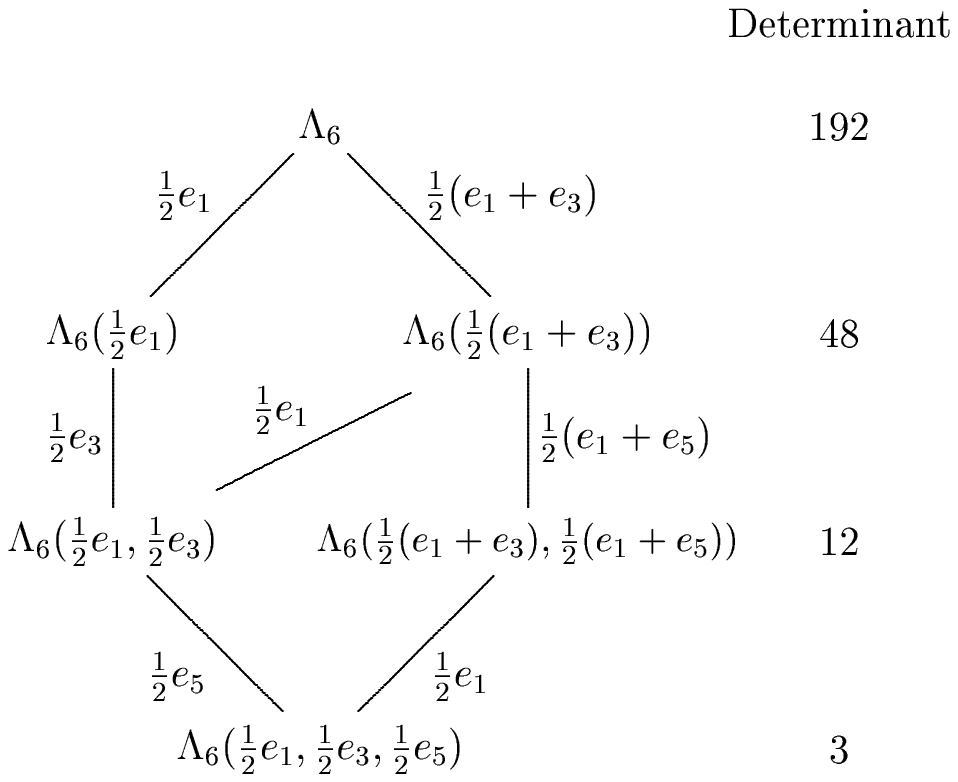}\vspace{-0.1in}
\caption{Lattice neighboring from $\Lambda_6$}\label{fig6}
\end{figure}

Here, a Gram matrix of $\Lambda_6$ is
\begin{equation*}
{}^t \! (e_1, \ldots, e_6) \cdot (e_1, \ldots, e_6) = 2
 {\small \left[ \begin{array}{cccccc} 2 & -1 & 0 & 0 & 0 & 0\\ -1 & 2 & -1 & 0 & 0 & 0\\ 0 & -1 & 2 & -1 & 0 & -1\\ 0 & 0 & -1 & 2 & -1 & 0\\ 0 & 0 & 0 & -1 & 2 & 0\\ 0 & 0 & -1 & 0 & 0 & 2 \end{array} \right]},
\end{equation*}
we have $\Lambda_6 (\frac{1}{2} (e_1 + e_3), \frac{1}{2} (e_1 + e_5)) \simeq \bm{A}_2 \perp \bm{D}_4$ and $\Lambda_6 (\frac{1}{2} e_1, \frac{1}{2} e_3, \frac{1}{2} e_5) \simeq \mathbb{Z}^4 \perp \bm{A}_2$.

We can show the following fact:
\begin{theorem}
The integral lattices which include $\Lambda_6$ are isometric to one of the above six lattices $($in Figure $\ref{fig6}$.$)$, namely $\Lambda_6$, $\Lambda_6 (\frac{1}{2} e_1)$, $\Lambda_6 (\frac{1}{2} (e_1 + e_3))$, $\Lambda_6 (\frac{1}{2} e_1, \frac{1}{2} e_3)$, $\Lambda_6 (\frac{1}{2} (e_1 + e_3), \frac{1}{2} (e_1 + e_5))$, and $\Lambda_6 (\frac{1}{2} e_1, \frac{1}{2} e_3, \frac{1}{2} e_5)$.
\end{theorem}

It is easy to classify the above lattices. We have the following:
\begin{proposition}
The $192$ classes of $\Lambda_6^{\sharp} / \Lambda_6$ consisit of the following classes:\\

\begin{center}
\begin{tabular}{ccc}
\hline
\# of class & norm of c.l. & order\\
\hline
$1$ & $0$ & $1$\\
$27 \times 2$ & $2/3$ & $6$\\
$36$ & $1$ & $2$\\
$36 \times 2$ & $5/3$ & $6$\\
$27$ & $2$ & $2$\\
$1 \times 2$ & $8/3$ & $3$\\
\hline
\end{tabular}
\end{center}\quad
\end{proposition}

Thus, to classify such lattices, we have only to consider the neighboring with the each vector of norm $1$ and $2$. We can easily check by computer search. In fact, we have $514$ integral lattices from the combinations of the above classes. The following table contains the number of integral lattices which are isometric to each lattice:\\

\begin{center}
\begin{tabular}{c|ccc|ccc|ccc|c}
\cline{1-2} \cline{4-5} \cline{7-8} \cline{10-11}
$\Lambda_6$ & $1$ & \qquad & $\Lambda_6 (\frac{1}{2} e_1)$ & $36$ & \qquad & $\Lambda_6 (\frac{1}{2} e_1, \frac{1}{2} e_3)$ & $270$ & \qquad & $\Lambda_6 (\frac{1}{2} e_1, \frac{1}{2} e_3, \frac{1}{2} e_5)$ & $135$\\
  &&& $\Lambda_6 (\frac{1}{2} (e_1 + e_3))$ & $27$ && {\small $\Lambda_6 (\frac{1}{2} (e_1 + e_3), \frac{1}{2} (e_1 + e_5))$} & $45$ &&&\\
\cline{1-2} \cline{4-5} \cline{7-8} \cline{10-11}
\end{tabular}
\end{center}\quad

\subsection{Properties of lattices}

The theta series of each lattice have the following form:
\begin{align*}
\Theta_{\Lambda_6} &= 1 + 72 \, q^4 + 270 \, q^8 + 720 \, q^{12} + \cdots\\
\Theta_{\Lambda_6 (\frac{1}{2} (e_1 + e_3))} &= 1 + 10 \, q^2 + 104 \, q^4 + 82 \, q^6 + 430 \, q^8 + 192 \, q^{10} + 1040 \, q^{12} + \cdots\\
\Theta_{\Lambda_6 (\frac{1}{2} (e_1 + e_3), \frac{1}{2} (e_1 + e_5))} &= 1 + 30 \, q^2 + 168 \, q^4 + 246 \, q^6 + 750 \, q^8 + 576 \, q^{10} + 1680 \, q^{12} + \cdots\\
\quad\\
\Theta_{\Lambda_6 (\frac{1}{2} e_1)} &= 1 + 2 \, q + 20 \, q^3 + 72 \, q^4 + 60 \, q^5 + 100 \, q^7\\
 &+ 270 \, q^8 + 182 \, q^9 + 300 \, q^{11} + 720 \, q^{12} + \cdots\\
\Theta_{\Lambda_6 (\frac{1}{2} e_1, \frac{1}{2} e_3)} &= 1 + 4 \, q + 10 \, q^2 + 40 \, q^3 + 104 \, q^4 + 120 \, q^5 + 82 \, q^6 + 200 \, q^7\\
 &+ 430 \, q^8 + 364 \, q^9 + 192 \, q^{10} + 600 \, q^{11} + 1040 \, q^{12} + \cdots\\
\Theta_{\Lambda_6 (\frac{1}{2} e_1, \frac{1}{2} e_3, \frac{1}{2} e_5)} &= 1 + 8 \, q + 30 \, q^2 + 80 \, q^3 + 168 \, q^4 + 240 \, q^5 + 246 \, q^6 + 400 \, q^7\\
 &+ 750 \, q^8 + 728 \, q^9 + 576 \, q^{10} + 1200 \, q^{11} + 1680 \, q^{12} + \cdots
\end{align*}\quad

The following table contains the $(d, n, s, t)$-configuration of each shell of norm $m$ of the lattice:
\begin{center}
\begin{tabular}{r}
\quad\\
\hline
$m$\\
\hline
$2$\\
$4$\\
$6$\\
$8$\\
$10$\\
$12$\\
\hline
\end{tabular}
\quad
\begin{tabular}{cccc}
\multicolumn{4}{c}{$\Lambda_6$}\\
\hline
$d$ & $n$ & $s$ & $t$\\
\hline
\ \\
$6$ & $72$ & $4$ & $5$\\
\ \\
$6$ & $270$ & $8$ & $5$\\
\ \\
$6$ & $720$ & $12$ & $5$\\
\hline
\end{tabular}
\qquad
\begin{tabular}{cccc}
\multicolumn{4}{c}{$\Lambda_6 (\frac{1}{2} (e_1 + e_3))$}\\
\hline
$d$ & $n$ & $s$ & $t$\\
\hline
$5$ & $10$ & $2$ & $3$\\
$6$ & $104$ & $8$ & $1$\\
$6$ & $82$ & $6$ & $1$\\
$6$ & $430$ & $16$ & $1$\\
$6$ & $192$ & $10$ & $1$\\
$6$ & $1040$ & $24$ & $1$\\
\hline
\end{tabular}
\qquad
\begin{tabular}{cccc}
\multicolumn{4}{c}{\tiny $\Lambda_6 (\frac{1}{2} (e_1 + e_3), \frac{1}{2} (e_1 + e_5))$}\\
\hline
\, $d$ \, & \, $n$ \, & \, $s$ \, & \, $t$ \, \\
\hline
$6$ & $30$ & $4$ & $1$\\
$6$ & $168$ & $8$ & $1$\\
$6$ & $246$ & $12$ & $1$\\
$6$ & $750$ & $16$ & $1$\\
$6$ & $576$ & $20$ & $1$\\
$6$ & $1680$ & $24$ & $1$\\
\hline
\end{tabular}
\end{center}\quad

\begin{center}
\begin{tabular}{r}
\quad\\
\hline
$m$\\
\hline
$1$\\
$2$\\
$3$\\
$4$\\
$5$\\
$6$\\
$7$\\
$8$\\
$9$\\
$10$\\
$11$\\
$12$\\
\hline
\end{tabular}
\qquad
\begin{tabular}{cccc}
\multicolumn{4}{c}{$\Lambda_6 (\frac{1}{2} e_1)$}\\
\hline
$d$ & $n$ & $s$ & $t$\\
\hline
$1$ & $2$ & $1$ & $-$\\
\ \\
$5$ & $20$ & $3$ & $3$\\
$6$ & $72$ & $4$ & $5$\\
$6$ & $60$ & $5$ & $1$\\
\ \\
$6$ & $100$ & $7$ & $1$\\
$6$ & $270$ & $8$ & $5$\\
$6$ & $182$ & $9$ & $1$\\
\ \\
$6$ & $300$ & $11$ & $1$\\
$6$ & $720$ & $12$ & $5$\\
\hline
\end{tabular}
\qquad
\begin{tabular}{cccc}
\multicolumn{4}{c}{$\Lambda_6 (\frac{1}{2} e_1, \frac{1}{2} e_3)$}\\
\hline
$d$ & $n$ & $s$ & $t$\\
\hline
$2$ & $4$ & $2$ & $3$\\
$5$ & $10$ & $2$ & $3$\\
$6$ & $40$ & $6$ & $1$\\
$6$ & $104$ & $8$ & $1$\\
$6$ & $120$ & $10$ & $1$\\
$6$ & $82$ & $6$ & $1$\\
$6$ & $200$ & $14$ & $1$\\
$6$ & $430$ & $16$ & $1$\\
$6$ & $364$ & $18$ & $1$\\
$6$ & $192$ & $10$ & $1$\\
$6$ & $600$ & $22$ & $1$\\
$6$ & $1040$ & $24$ & $1$\\
\hline
\end{tabular}
\qquad
\begin{tabular}{cccc}
\multicolumn{4}{c}{\small $\Lambda_6 (\frac{1}{2} e_1, \frac{1}{2} e_3, \frac{1}{2} e_5)$}\\
\hline
$d$ & $n$ & $s$ & $t$\\
\hline
$4$ & $8$ & $2$ & $3$\\
$6$ & $30$ & $4$ & $1$\\
$6$ & $80$ & $6$ & $1$\\
$6$ & $168$ & $8$ & $1$\\
$6$ & $240$ & $10$ & $1$\\
$6$ & $246$ & $12$ & $1$\\
$6$ & $400$ & $14$ & $1$\\
$6$ & $750$ & $16$ & $1$\\
$6$ & $728$ & $18$ & $1$\\
$6$ & $576$ & $20$ & $1$\\
$6$ & $1200$ & $22$ & $1$\\
$6$ & $1680$ & $24$ & $1$\\
\hline
\end{tabular}
\end{center}

\newpage

\section{The laminated lattice $\Lambda_7 \simeq \sqrt{2} \, \bm{E}_7$}

We have the following lattice neighboring from the laminated lattice $\Lambda_7$:\vspace{-0.1in}

\begin{figure}[h]
\includegraphics[width=5.2in]{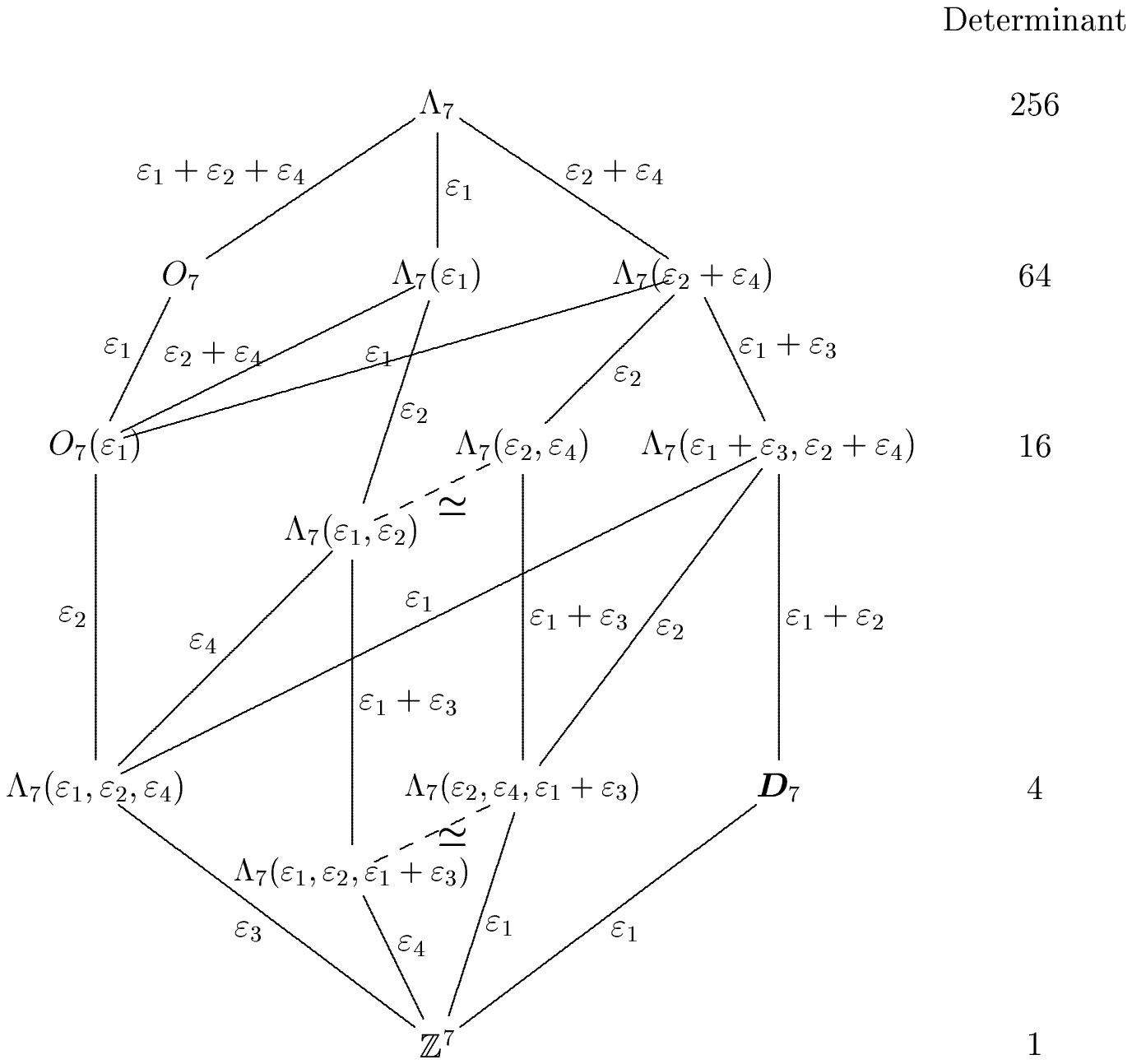}\vspace{-0.1in}
\caption{Lattice neighboring from $\Lambda_7$}\label{fig7}
\end{figure}
Note that $\Lambda_7 (\varepsilon_1, \varepsilon_2) \simeq \Lambda_7 (\varepsilon_2, \varepsilon_4)$ and $\Lambda_7 (\varepsilon_1, \varepsilon_2, \varepsilon_1 + \varepsilon_3) \simeq \Lambda_7 (\varepsilon_2, \varepsilon_4, \varepsilon_1 + \varepsilon_3)$.\\

Here, for the above neighboring, we have the following generator matrix of $\Lambda_7$:
{\small \begin{equation*}
\left(\begin{array}{ccccccc}
2 & 0 & 0 & 0 & 0 & 0 & 0\\
0 & 2 & 0 & 0 & 0 & 0 & 0\\
0 & 0 & 2 & 0 & 0 & 0 & 0\\
0 & 0 & 0 & 2 & 0 & 0 & 0\\
1 & 0 & 1 & 1 & 1 & 0 & 0\\
0 & 1 & 0 & 1 & 1 & 1 & 0\\
0 & 0 & 1 & 0 & 1 & 1 & 1
\end{array}\right)
\end{equation*}}
we have $\Lambda_7 (\varepsilon_1, \varepsilon_2, \varepsilon_4) \simeq \mathbb{Z}^3 \perp \bm{D}_4$ and $\Lambda_7 (\varepsilon_1, \varepsilon_2, \varepsilon_1 + \varepsilon_3) \simeq \mathbb{Z}^4 \perp \bm{A}_3$.

We can show the following fact:
\begin{theorem}
The integral lattices which include $\Lambda_7$ are isometric to one of the above eleven lattices $($in Figure $\ref{fig7}$.$)$, namely $\Lambda_7$, $O_7$, $\Lambda_7 (\varepsilon_1)$, $\Lambda_7 (\varepsilon_2 + \varepsilon_4)$, $O_7 (\varepsilon_1)$, $\Lambda_7 (\varepsilon_1, \varepsilon_2)$, $\Lambda_7 (\varepsilon_1 + \varepsilon_3, \varepsilon_2 + \varepsilon_4)$, $\Lambda_7 (\varepsilon_1, \varepsilon_2, \varepsilon_4)$, $\Lambda_7 (\varepsilon_1, \varepsilon_2, \varepsilon_1 + \varepsilon_3)$, $\bm{D}_7$, and $\mathbb{Z}^7$.
\end{theorem}

\newpage

It is easy to classify the above lattices. We have the following:
\begin{proposition}
The $256$ classes of $\Lambda_7^{\sharp} / \Lambda_7$ consisit of the following classes:\\

\begin{center}
\begin{tabular}{ccc}
\hline
\# of class & norm of c.l. & order\\
\hline
$1$ & $0$ & $1$\\
$28 \times 2$ & $3/4$ & $4$\\
$63$ & $1$ & $2$\\
$36 \times 2$ & $7/4$ & $4$\\
$63$ & $2$ & $2$\\
$1$ & $3$ & $2$\\
\hline
\end{tabular}
\end{center}\quad
\end{proposition}

Thus, to classify such lattices, we have only to consider the neighboring with the each vector of norm $1$, $2$, and $3$. We can easily check by computer search. In fact, we have $2981$ integral lattices from the combinations of the above classes. The following table contains the number of integral lattices which are isometric to each lattice:\\

\begin{center}
\begin{tabular}{c|ccc|ccc|ccc|ccc|c}
\cline{1-2} \cline{4-5} \cline{7-8} \cline{10-11} \cline{13-14}
$\Lambda_7$ & $1$ & \qquad & $O_7$ & $1$ & \qquad & $O_7 (\varepsilon_1)$ & $63$ & \qquad & $\Lambda_7 (\varepsilon_1, \varepsilon_2, \varepsilon_4)$ & $315$ & \qquad & $\mathbb{Z}^7$ & $135$\\
 &&& $\Lambda_7 (\varepsilon_1)$ & $63$ && $\Lambda_7 (\varepsilon_1, \varepsilon_2)$ & $945$ && $\Lambda_7 (\varepsilon_1, \varepsilon_2, \varepsilon_1 + \varepsilon_3)$ & $945$ &&&\\
  &&& $\Lambda_7 (\varepsilon_2 + \varepsilon_4)$ & $63$ && {\small $\Lambda_7 (\varepsilon_1 + \varepsilon_3, \varepsilon_2 + \varepsilon_4)$} & $315$ && $\bm{D}_7$ & $135$ &&&\\
\cline{1-2} \cline{4-5} \cline{7-8} \cline{10-11} \cline{13-14}
\end{tabular}
\end{center}\quad

\subsection{Properties of lattices}

The theta series of each lattice have the following form:
\begin{align*}
\Theta_{\Lambda_7} &= 1 + 126 \, q^4 + 756 \, q^8 + 2072 \, q^{12} + \cdots\\
\quad\\
\Theta_{\Lambda_7 (\varepsilon_2 + \varepsilon_4)} &= 1 + 12 \, q^2 + 190 \, q^4 + 184 \, q^6 + 1140 \, q^8 + 632 \, q^{10} + 3096 \, q^{12} + \cdots\\
\Theta_{\Lambda_7 (\varepsilon_1 + \varepsilon_3, \varepsilon_2 + \varepsilon_4)} &= 1 + 36 \, q^2 + 318 \, q^4 + 552 \, q^6 + 1908 \, q^8 + 1896 \, q^{10} + 5144 \, q^{12} + \cdots\\
\Theta_{\bm{D}_7} &= 1 + 84 \, q^2 + 574 \, q^4 + 1288 \, q^6 + 3444 \, q^8 + 4424 \, q^{10} + 9240 \, q^{12} + \cdots\\
\quad\\
\Theta_{O_7} &= 1 + 56 \, q^3 + 126 \, q^4 + 576 \, q^7 + 756 \, q^8 + 1512 \, q^{11} + 2072 \, q^{12} + \cdots\\
\Theta_{O_7 (\varepsilon_1)} &= 1 + 2 \, q + 12 \, q^2 + 88 \, q^3 + 190 \, q^4 + 120 \, q^5 + 184 \, q^6 + 832 \, q^7\\
 &+ 1140 \, q^8 + 506 \, q^9 + 632 \, q^{10} + 2376 \, q^{11} + 3096 \, q^{12} + \cdots\\
\Theta_{\Lambda_7 (\varepsilon_1, \varepsilon_2, \varepsilon_4)} &= 1 + 6 \, q + 36 \, q^2 + 152 \, q^3 + 318 \, q^4 + 360 \, q^5 + 552 \, q^6 + 1344 \, q^7\\
 &+ 1908 \, q^8 + 1518 \, q^9 + 1896 \, q^{10} + 4104 \, q^{11} + 5144 \, q^{12} + \cdots\\
\quad\\
\Theta_{\Lambda_7 (\varepsilon_1)} &= 1 + 2 \, q + 32 \, q^3 + 126 \, q^4 + 120 \, q^5 + 256 \, q^7\\
 &+ 756 \, q^8 + 506 \, q^9 + 864 \, q^{11} + 2072 \, q^{12} + \cdots\\
\Theta_{\Lambda_7 (\varepsilon_1, \varepsilon_2)} &= 1 + 4 \, q + 12 \, q^2 + 64 \, q^3 + 190 \, q^4 + 240 \, q^5 + 184 \, q^6 + 512 \, q^7\\
 &+ 1140 \, q^8 + 1012 \, q^9 + 632 \, q^{10} + 1728 \, q^{11} + 3096 \, q^{12} + \cdots\\
\Theta_{\Lambda_7 (\varepsilon_1, \varepsilon_2, \varepsilon_1 + \varepsilon_3)} &= 1 + 8 \, q + 36 \, q^2 + 128 \, q^3 + 318 \, q^4 + 480 \, q^5 + 552 \, q^6 + 1024 \, q^7\\
 &+ 1908 \, q^8 + 2024 \, q^9 + 1896 \, q^{10} + 3456 \, q^{11} + 5144 \, q^{12} + \cdots\\
\quad\\
\Theta_{\mathbb{Z}^7} &= 1 + 14 \, q + 84 \, q^2 + 280 \, q^3 + 574 \, q^4 + 840 \, q^5 + 1288 \, q^6 + 2368 \, q^7\\
 &+ 3444 \, q^8 + 3542 \, q^9 + 4424 \, q^{10} + 7560 \, q^{11} + 9240 \, q^{12} + \cdots\\
\end{align*}

\newpage

The following table contains the $(d, n, s, t)$-configuration of each shell of norm $m$ of the lattice:
\begin{center}
\begin{tabular}{r}
\quad\\
\hline
$m$\\
\hline
$2$\\
$4$\\
$6$\\
$8$\\
$10$\\
$12$\\
\hline
\end{tabular}
\quad
\begin{tabular}{cccc}
\multicolumn{4}{c}{$\Lambda_7$}\\
\hline
$d$ & $n$ & $s$ & $t$\\
\hline
\ \\
$7$ & $126$ & $4$ & $5$\\
\ \\
$7$ & $756$ & $8$ & $5$\\
\ \\
$7$ & $2072$ & $12$ & $5$\\
\hline
\end{tabular}
\qquad
\begin{tabular}{cccc}
\multicolumn{4}{c}{$\Lambda_7 (\varepsilon_2 + \varepsilon_4)$}\\
\hline
$d$ & $n$ & $s$ & $t$\\
\hline
$6$ & $12$ & $2$ & $3$\\
$7$ & $190$ & $8$ & $1$\\
$7$ & $184$ & $6$ & $1$\\
$7$ & $1140$ & $16$ & $1$\\
$7$ & $632$ & $10$ & $1$\\
$7$ & $3096$ & $24$ & $1$\\
\hline
\end{tabular}
\qquad
\begin{tabular}{cccc}
\multicolumn{4}{c}{\small $\Lambda_7 (\varepsilon_1 + \varepsilon_3, \varepsilon_2 + \varepsilon_4)$}\\
\hline
$d$ & $n$ & $s$ & $t$\\
\hline
$7$ & $36$ & $4$ & $1$\\
$7$ & $318$ & $8$ & $1$\\
$7$ & $552$ & $12$ & $1$\\
$7$ & $1908$ & $16$ & $1$\\
$7$ & $1896$ & $20$ & $1$\\
$7$ & $5144$ & $24$ & $1$\\
\hline
\end{tabular}
\qquad
\begin{tabular}{cccc}
\multicolumn{4}{c}{$\bm{D}_7$}\\
\hline
$d$ & $n$ & $s$ & $t$\\
\hline
$7$ & $84$ & $4$ & $3$\\
$7$ & $574$ & $8$ & $3$\\
$7$ & $1288$ & $12$ & $3$\\
$7$ & $3444$ & $16$ & $3$\\
$7$ & $4424$ & $20$ & $3$\\
$7$ & $9240$ & $24$ & $5$\\
\hline
\end{tabular}
\end{center}\quad

\begin{center}
\begin{tabular}{r}
\quad\\
\hline
$m$\\
\hline
$1$\\
$2$\\
$3$\\
$4$\\
$5$\\
$6$\\
$7$\\
$8$\\
$9$\\
$10$\\
$11$\\
$12$\\
\hline
\end{tabular}
\qquad
\begin{tabular}{cccc}
\multicolumn{4}{c}{$O_7$}\\
\hline
$d$ & $n$ & $s$ & $t$\\
\hline
\ \\ \ \\
$7$ & $56$ & $3$ & $5$\\
$7$ & $126$ & $4$ & $5$\\
\ \\ \ \\
$7$ & $576$ & $7$ & $5$\\
$7$ & $756$ & $8$ & $5$\\
\ \\ \ \\
$7$ & $1512$ & $11$ & $5$\\
$7$ & $2072$ & $12$ & $5$\\
\hline
\end{tabular}
\qquad
\begin{tabular}{cccc}
\multicolumn{4}{c}{$O_7 (\varepsilon_1)$}\\
\hline
$d$ & $n$ & $s$ & $t$\\
\hline
$1$ & $2$ & $1$ & $-$\\
$6$ & $12$ & $2$ & $3$\\
$7$ & $88$ & $6$ & $1$\\
$7$ & $190$ & $8$ & $1$\\
$7$ & $120$ & $5$ & $1$\\
$7$ & $184$ & $6$ & $1$\\
$7$ & $832$ & $14$ & $3$\\
$7$ & $1140$ & $16$ & $1$\\
$7$ & $506$ & $9$ & $1$\\
$7$ & $632$ & $10$ & $1$\\
$7$ & $2376$ & $22$ & $1$\\
$7$ & $3096$ & $24$ & $1$\\
\hline
\end{tabular}
\qquad
\begin{tabular}{cccc}
\multicolumn{4}{c}{\small $\Lambda_7 (\varepsilon_1, \varepsilon_2, \varepsilon_4)$}\\
\hline
$d$ & $n$ & $s$ & $t$\\
\hline
$3$ & $6$ & $2$ & $3$\\
$7$ & $36$ & $4$ & $1$\\
$7$ & $152$ & $6$ & $1$\\
$7$ & $318$ & $8$ & $1$\\
$7$ & $360$ & $10$ & $1$\\
$7$ & $552$ & $12$ & $1$\\
$7$ & $1344$ & $14$ & $3$\\
$7$ & $1908$ & $16$ & $1$\\
$7$ & $1518$ & $18$ & $1$\\
$7$ & $1896$ & $20$ & $1$\\
$7$ & $4104$ & $22$ & $1$\\
$7$ & $5144$ & $24$ & $1$\\
\hline
\end{tabular}
\end{center}\quad

\begin{center}
\begin{tabular}{r}
\quad\\
\hline
$m$\\
\hline
$1$\\
$2$\\
$3$\\
$4$\\
$5$\\
$6$\\
$7$\\
$8$\\
$9$\\
$10$\\
$11$\\
$12$\\
\hline
\end{tabular}
\qquad
\begin{tabular}{cccc}
\multicolumn{4}{c}{$\Lambda_7 (\varepsilon_1)$}\\
\hline
$d$ & $n$ & $s$ & $t$\\
\hline
$1$ & $2$ & $1$ & $-$\\
\ \\
$6$ & $32$ & $3$ & $3$\\
$7$ & $126$ & $4$ & $5$\\
$7$ & $120$ & $5$ & $1$\\
\ \\
$7$ & $256$ & $7$ & $3$\\
$7$ & $756$ & $8$ & $5$\\
$7$ & $506$ & $9$ & $1$\\
\ \\
$7$ & $864$ & $11$ & $1$\\
$7$ & $2072$ & $12$ & $5$\\
\hline
\end{tabular}
\qquad
\begin{tabular}{cccc}
\multicolumn{4}{c}{$\Lambda_7 (\varepsilon_1, \varepsilon_2)$}\\
\hline
$d$ & $n$ & $s$ & $t$\\
\hline
$2$ & $4$ & $2$ & $3$\\
$6$ & $12$ & $2$ & $3$\\
$7$ & $64$ & $6$ & $1$\\
$7$ & $190$ & $8$ & $1$\\
$7$ & $240$ & $10$ & $1$\\
$7$ & $184$ & $6$ & $1$\\
$7$ & $512$ & $14$ & $3$\\
$7$ & $1140$ & $16$ & $1$\\
$7$ & $1012$ & $18$ & $1$\\
$7$ & $632$ & $10$ & $1$\\
$7$ & $1728$ & $22$ & $1$\\
$7$ & $3096$ & $24$ & $1$\\
\hline
\end{tabular}
\qquad
\begin{tabular}{cccc}
\multicolumn{4}{c}{\small $\Lambda_7 (\varepsilon_1, \varepsilon_2, \varepsilon_1 + \varepsilon_3)$}\\
\hline
$d$ & $n$ & $s$ & $t$\\
\hline
$4$ & $8$ & $2$ & $3$\\
$7$ & $36$ & $4$ & $1$\\
$7$ & $128$ & $6$ & $1$\\
$7$ & $318$ & $8$ & $1$\\
$7$ & $480$ & $10$ & $1$\\
$7$ & $552$ & $12$ & $1$\\
$7$ & $1024$ & $14$ & $3$\\
$7$ & $1908$ & $16$ & $1$\\
$7$ & $2024$ & $18$ & $1$\\
$7$ & $1896$ & $20$ & $1$\\
$7$ & $3456$ & $22$ & $1$\\
$7$ & $5144$ & $24$ & $1$\\
\hline
\end{tabular}
\qquad
\begin{tabular}{cccc}
\multicolumn{4}{c}{$\mathbb{Z}^7$}\\
\hline
$d$ & $n$ & $s$ & $t$\\
\hline
$7$ & $14$ & $2$ & $3$\\
$7$ & $84$ & $4$ & $3$\\
$7$ & $280$ & $6$ & $5$\\
$7$ & $574$ & $8$ & $3$\\
$7$ & $840$ & $10$ & $3$\\
$7$ & $1288$ & $12$ & $3$\\
$7$ & $2368$ & $14$ & $3$\\
$7$ & $3444$ & $16$ & $3$\\
$7$ & $3542$ & $18$ & $3$\\
$7$ & $4424$ & $20$ & $3$\\
$7$ & $7560$ & $22$ & $5$\\
$7$ & $9240$ & $24$ & $5$\\
\hline
\end{tabular}
\end{center}

\newpage

\section{The laminated lattice $\Lambda_8 \simeq \sqrt{2} \, \bm{E}_8$}

We have the following lattice neighboring from the laminated lattice $\Lambda_8$:\vspace{-0.1in}

\begin{figure}[h]
\includegraphics[width=2.4in]{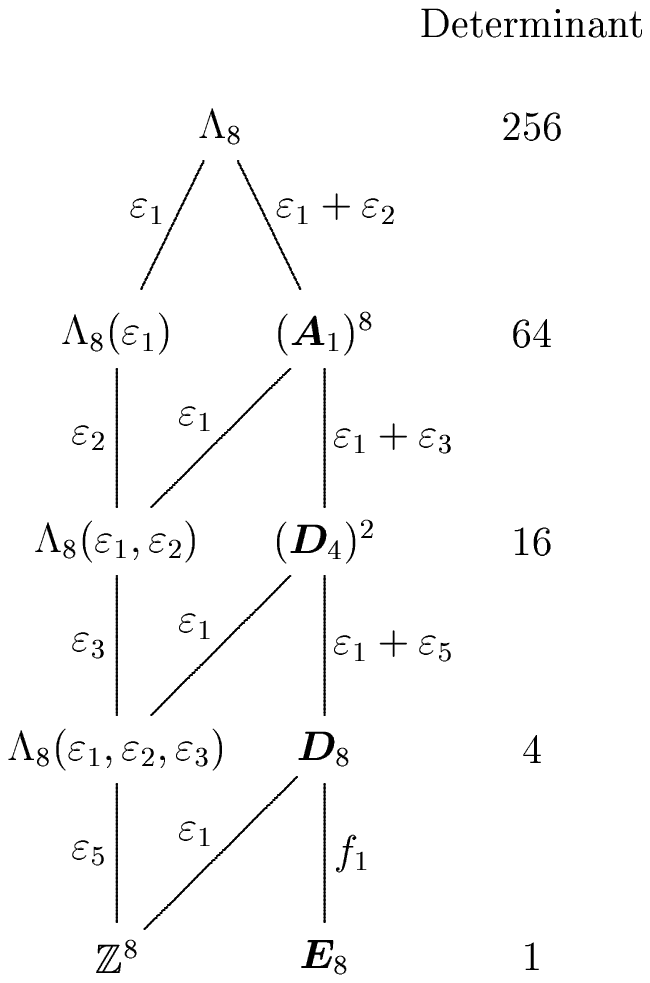}\vspace{-0.1in}
\caption{Lattice neighboring from $\Lambda_8$}\label{fig8}
\end{figure}

Here, for the above neighboring, we have the following generator matrix of $\Lambda_8$:
{\small \begin{equation*}
\left(\begin{array}{cccccccc}
2 & 0 & 0 & 0 & 0 & 0 & 0 & 0\\
0 & 2 & 0 & 0 & 0 & 0 & 0 & 0\\
0 & 0 & 2 & 0 & 0 & 0 & 0 & 0\\
1 & 1 & 1 & 1 & 0 & 0 & 0 & 0\\
0 & 0 & 0 & 0 & 2 & 0 & 0 & 0\\
1 & 1 & 0 & 0 & 1 & 1 & 0 & 0\\
1 & 0 & 1 & 0 & 1 & 0 & 1 & 0\\
0 & 1 & 0 & 1 & 0 & 1 & 0 & 1
\end{array}\right)
\end{equation*}}
We also have $\Lambda_8 (\varepsilon_1) \simeq \mathbb{Z} \perp O_7$ and $\Lambda_8 (\varepsilon_1, \varepsilon_2, \varepsilon_3) \simeq \mathbb{Z}^4 \perp \bm{D}_4$.

We can show the following fact:
\begin{theorem}
The integral lattices which include $\Lambda_8$ are isometric to one of the above eleven lattices $($in Figure $\ref{fig8}$.$)$, namely $\Lambda_8$, $O_8$, $\Lambda_8 (\varepsilon_1)$, $\Lambda_8 (\varepsilon_1, \varepsilon_2)$, $\Lambda_8 (\varepsilon_1, \varepsilon_2, \varepsilon_3)$, $\mathbb{Z}^8$, $(\bm{A}_1)^8$, $(\bm{D}_4)^2$, $\bm{D}_8$, and $\mathbb{Z}^8$.
\end{theorem}

It is easy to classify the above lattices. We have the following:
\begin{proposition}
The $256$ classes of $\Lambda_8^{\sharp} / \Lambda_8$ consisit of the following classes:\\

\begin{center}
\begin{tabular}{ccc}
\hline
\# of class & norm of c.l. & order\\
\hline
$1$ & $0$ & $1$\\
$120$ & $1$ & $2$\\
$135$ & $2$ & $2$\\
\hline
\end{tabular}
\end{center}\quad
\end{proposition}

Thus, to classify such lattices, we have only to consider the neighboring with the each vector of norm $1$, and $2$. We can easily check by computer search. In fact, we have $19381$ integral lattices from the combinations of the above classes. The following table contains the number of integral lattices which are isometric to each lattice:\\

\begin{center}
\begin{tabular}{c|ccc|ccc|ccc|ccc|c}
\cline{1-2} \cline{4-5} \cline{7-8} \cline{10-11} \cline{13-14}
$\Lambda_8$ & $1$ & \qquad & $\Lambda_8 (\varepsilon_1)$ & $120$ & \qquad & $\Lambda_8 (\varepsilon_1, \varepsilon_2)$ & $3780$ & \qquad & $\Lambda_8 (\varepsilon_1, \varepsilon_2, \varepsilon_3)$ & $9450$ & \qquad & $\mathbb{Z}^8$ & $2025$\\
  &&& $(\bm{A}_1)^8$ & $135$ && $(\bm{D}_4)^2$ & $1575$ && $\bm{D}_8$ & $2025$ && $\bm{E}_8$ & $270$\\
\cline{1-2} \cline{4-5} \cline{7-8} \cline{10-11} \cline{13-14}
\end{tabular}
\end{center}\quad

\subsection{Properties of lattices}

The theta series of each lattice have the following form:
\begin{align*}
\Theta_{\Lambda_8} &= 1 + 240 \, q^4 + 2160 \, q^8 + 6720 \, q^{12} + \cdots\\
\Theta_{(\bm{A}_1)^8} &= 1 + 16 \, q^2 + 368 \, q^4 + 448 \, q^6 + 3184 \, q^8 + 2016 \, q^{10} + 10304 \, q^{12} + \cdots\\
\Theta_{(\bm{D}_4)^2} &= 1 + 48 \, q^2 + 624 \, q^4 + 1344 \, q^6 + 5232 \, q^8 + 6048 \, q^{10} + 17472 \, q^{12} + \cdots\\
\Theta_{\bm{D}_8} &= 1 + 112 \, q^2 + 1136 \, q^4 + 3136 \, q^6 + 9328 \, q^8 + 14112 \, q^{10} + 31808 \, q^{12} + \cdots\\
\Theta_{\bm{E}_8} &= 1 + 240 \, q^2 + 2160 \, q^4 + 6720 \, q^6 + 17520 \, q^8 + 30240 \, q^{10} + 60480 \, q^{12} + \cdots\\
\quad\\
\Theta_{\Lambda_8 (\varepsilon_1)} &= 1 + 2 \, q + 56 \, q^3 + 240 \, q^4 + 252 \, q^5 + 688 \, q^7\\
 &+ 2160 \, q^8 + 1514 \, q^9 + 2664 \, q^{11} + 6720 \, q^{12} + \cdots\\
\Theta_{\Lambda_8 (\varepsilon_1, \varepsilon_2)} &= 1 + 4 \, q + 16 \, q^2 + 112 \, q^3 + 368 \, q^4 + 504 \, q^5 + 448 \, q^6 + 1376 \, q^7\\
 &+ 3184 \, q^8 + 3028 \, q^9 + 2016 \, q^{10} + 5328 \, q^{11} + 10304 \, q^{12} + \cdots\\
\Theta_{\Lambda_8 (\varepsilon_1, \varepsilon_2, \varepsilon_3)} &= 1 + 8 \, q + 48 \, q^2 + 224 \, q^3 + 624 \, q^4 + 1008 \, q^5 + 1344 \, q^6 + 2752 \, q^7\\
 &+ 5232 \, q^8 + 6056 \, q^9 + 6048 \, q^{10} + 10656 \, q^{11} + 17472 \, q^{12} + \cdots\\
\Theta_{\mathbb{Z}^8} &= 1 + 16 \, q + 112 \, q^2 + 448 \, q^3 + 1136 \, q^4 + 2016 \, q^5 + 3136 \, q^6 + 5504 \, q^7\\
 &+ 9328 \, q^8 + 12112 \, q^9 + 14112 \, q^{10} + 21312 \, q^{11} + 31808 \, q^{12} + \cdots
\end{align*}\quad

The following table contains the $(d, n, s, t)$-configuration of each shell of norm $m$ of the lattice:
\begin{center}
\begin{tabular}{r}
\quad\\
\hline
$m$\\
\hline
$4$\\
$8$\\
$12$\\
\hline
\end{tabular}
\quad
\begin{tabular}{cccc}
\multicolumn{4}{c}{$\Lambda_8$}\\
\hline
$d$ & $n$ & $s$ & $t$\\
\hline
$8$ & $240$ & $4$ & $7$\\
$8$ & $2160$ & $8$ & $7$\\
$8$ & $6720$ & $12$ & $7$\\
\hline
\end{tabular}
\end{center}\quad

\begin{center}
\begin{tabular}{r}
\quad\\
\hline
$m$\\
\hline
$2$\\
$4$\\
$6$\\
$8$\\
$10$\\
$12$\\
\hline
\end{tabular}
\quad
\begin{tabular}{cccc}
\multicolumn{4}{c}{$(\bm{A}_1)^8$}\\
\hline
$d$ & $n$ & $s$ & $t$\\
\hline
$8$ & $16$ & $2$ & $3$\\
$8$ & $368$ & $8$ & $3$\\
$8$ & $448$ & $6$ & $3$\\
$8$ & $3184$ & $16$ & $3$\\
$8$ & $2016$ & $10$ & $3$\\
$8$ & $10304$ & $24$ & $3$\\
\hline
\end{tabular}
\qquad
\begin{tabular}{cccc}
\multicolumn{4}{c}{$(\bm{D}_4)^2$}\\
\hline
$d$ & $n$ & $s$ & $t$\\
\hline
$8$ & $48$ & $4$ & $3$\\
$8$ & $624$ & $8$ & $3$\\
$8$ & $1344$ & $12$ & $3$\\
$8$ & $5234$ & $16$ & $3$\\
$8$ & $6048$ & $20$ & $3$\\
$8$ & $17472$ & $24$ & $3$\\
\hline
\end{tabular}
\qquad
\begin{tabular}{cccc}
\multicolumn{4}{c}{$\bm{D}_8$}\\
\hline
$d$ & $n$ & $s$ & $t$\\
\hline
$8$ & $112$ & $4$ & $3$\\
$8$ & $1136$ & $8$ & $3$\\
$8$ & $3136$ & $12$ & $3$\\
$8$ & $9328$ & $16$ & $3$\\
$8$ & $14112$ & $20$ & $3$\\
$8$ & $31808$ & $24$ & $3$\\
\hline
\end{tabular}
\qquad
\begin{tabular}{cccc}
\multicolumn{4}{c}{$\bm{E}_8$}\\
\hline
$d$ & $n$ & $s$ & $t$\\
\hline
$8$ & $240$ & $4$ & $7$\\
$8$ & $2160$ & $8$ & $7$\\
$8$ & $6720$ & $12$ & $7$\\
$8$ & $17520$ & $16$ & $7$\\
$8$ & $30240$ & $20$ & $7$\\
$8$ & $60480$ & $24$ & $7$\\
\hline
\end{tabular}
\end{center}\quad

\begin{center}
\begin{tabular}{r}
\quad\\
\hline
$m$\\
\hline
$1$\\
$2$\\
$3$\\
$4$\\
$5$\\
$6$\\
$7$\\
$8$\\
$9$\\
$10$\\
$11$\\
$12$\\
\hline
\end{tabular}
\qquad
\begin{tabular}{cccc}
\multicolumn{4}{c}{$\Lambda_8 (\varepsilon_1)$}\\
\hline
$d$ & $n$ & $s$ & $t$\\
\hline
$1$ & $2$ & $1$ & $-$\\
\ \\
$7$ & $56$ & $3$ & $5$\\
$8$ & $240$ & $4$ & $7$\\
$8$ & $252$ & $5$ & $1$\\
\ \\
$8$ & $688$ & $7$ & $1$\\
$8$ & $2160$ & $8$ & $7$\\
$8$ & $1514$ & $9$ & $1$\\
\ \\
$8$ & $2664$ & $11$ & $1$\\
$8$ & $6720$ & $12$ & $7$\\
\hline
\end{tabular}
\qquad
\begin{tabular}{cccc}
\multicolumn{4}{c}{$\Lambda_8 (\varepsilon_1, \varepsilon_2)$}\\
\hline
$d$ & $n$ & $s$ & $t$\\
\hline
$2$ & $4$ & $2$ & $3$\\
$8$ & $16$ & $2$ & $3$\\
$8$ & $112$ & $6$ & $1$\\
$8$ & $368$ & $8$ & $3$\\
$8$ & $504$ & $10$ & $1$\\
$8$ & $448$ & $6$ & $3$\\
$8$ & $1376$ & $14$ & $1$\\
$8$ & $3184$ & $16$ & $3$\\
$8$ & $3028$ & $18$ & $1$\\
$8$ & $2016$ & $10$ & $3$\\
$8$ & $5328$ & $22$ & $1$\\
$8$ & $10304$ & $24$ & $3$\\
\hline
\end{tabular}
\qquad
\begin{tabular}{cccc}
\multicolumn{4}{c}{\small $\Lambda_8 (\varepsilon_1, \varepsilon_2, \varepsilon_3)$}\\
\hline
$d$ & $n$ & $s$ & $t$\\
\hline
$4$ & $8$ & $2$ & $3$\\
$8$ & $48$ & $4$ & $3$\\
$8$ & $224$ & $6$ & $1$\\
$8$ & $624$ & $8$ & $3$\\
$8$ & $1008$ & $10$ & $1$\\
$8$ & $1344$ & $12$ & $3$\\
$8$ & $2752$ & $14$ & $1$\\
$8$ & $5232$ & $16$ & $3$\\
$8$ & $6056$ & $18$ & $1$\\
$8$ & $6048$ & $20$ & $3$\\
$8$ & $10656$ & $22$ & $1$\\
$8$ & $17472$ & $24$ & $3$\\
\hline
\end{tabular}
\qquad
\begin{tabular}{cccc}
\multicolumn{4}{c}{$\mathbb{Z}^8$}\\
\hline
$d$ & $n$ & $s$ & $t$\\
\hline
$8$ & $16$ & $2$ & $3$\\
$8$ & $112$ & $4$ & $3$\\
$8$ & $448$ & $6$ & $3$\\
$8$ & $1136$ & $8$ & $3$\\
$8$ & $2016$ & $10$ & $3$\\
$8$ & $3136$ & $12$ & $3$\\
$8$ & $5504$ & $14$ & $3$\\
$8$ & $9328$ & $16$ & $3$\\
$8$ & $12112$ & $18$ & $3$\\
$8$ & $14112$ & $20$ & $3$\\
$8$ & $21312$ & $22$ & $3$\\
$8$ & $31808$ & $24$ & $3$\\
\hline
\end{tabular}
\end{center}

\end{document}